\documentclass[12pt]{amsart}
\usepackage{amsmath,amsthm}
\usepackage{amssymb}
\usepackage{euscript}
\DeclareMathAlphabet{\EuRm}{U}{eur}{m}{n}
\SetMathAlphabet{\EuRm}{bold}{U}{eur}{b}{n}
\begin{document}
%
%
\newcounter{numb}
\newcounter{num}
%
%
\swapnumbers
\newtheorem{thm}{Theorem}[section]
\newtheorem*{tha}{Theorem A}
\newtheorem*{thb}{Theorem B}
\newtheorem{lemma}[thm]{Lemma}
\newtheorem{prop}[thm]{Proposition}
\newtheorem{cor}[thm]{Corollary}
\theoremstyle{definition}
\newtheorem{defn}[thm]{Definition}
\newtheorem{example}[thm]{Example}
\newtheorem{notation}[thm]{Notation}
\newtheorem{summary}[thm]{Summary}
\newtheorem{fact}[thm]{Fact}
\theoremstyle{remark}
\newtheorem{remark}[thm]{Remark}
\newtheorem{assume}[thm]{Assumption}
\newtheorem{note}[thm]{Note}
\newtheorem{ack}[thm]{Acknowledgements}

\numberwithin{equation}{section}
%
%
\def\sect{\setcounter{thm}{0} \section}
%
%
\newcommand{\xra}[1]{\xrightarrow{#1}}
\newcommand{\xla}[1]{\xleftarrow{#1}}
\newcommand{\hra}{\hookrightarrow}
\newcommand{\adj}[2]{\substack{{#1}\\ \rightleftharpoons \\ {#2}}}
\newcommand{\hsp}{\hspace{10 mm}}
\newcommand{\hs}{\hspace{5 mm}}
\newcommand{\hsm}{\hspace{2 mm}}
\newcommand{\vs}{\vspace{7 mm}}
\newcommand{\vsm}{\vspace{2 mm}}
\newcommand{\rest}[1]{\lvert_{#1}}
\newcommand{\lra}[1]{\langle{#1}\rangle}
\newcommand{\llra}[1]{\langle\!\langle{#1}\rangle\!\rangle}
\newcommand{\DEF}{:=}
\newcommand{\EQUIV}{\Leftrightarrow}
\newcommand{\epic}{\to\hspace{-5 mm}\to}
\newcommand{\xepic}[1]{\xrightarrow{#1}\hspace{-5 mm}\to}
\newcommand{\hotimes}{\hat{\otimes}}
%
%
\newcommand{\Coker}{\operatorname{Coker}}
\newcommand{\colim}{\operatorname{colim}}
\newcommand{\diag}{\operatorname{diag}}
\newcommand{\Ext}{\operatorname{Ext}}
\newcommand{\Fib}{\operatorname{Fib}}
\newcommand{\holim}{\operatorname{holim}}
\newcommand{\Hom}{\operatorname{Hom}}
\newcommand{\Image}{\operatorname{Im}}
\newcommand{\Ker}{\operatorname{Ker}}
\newcommand{\sk}[1]{\operatorname{sk}_{#1}}
\newcommand{\tr}[1]{\operatorname{tr}_{#1}}
%
%
\newcommand{\Ab}{{\EuScript Ab}}
\newcommand{\Alg}[1]{{#1}\text{-}{\EuScript Alg}}
\newcommand{\BB}{{\mathcal B}}
\newcommand{\C}{{\mathcal C}}
\newcommand{\D}{{\mathcal D}}
\newcommand{\F}{{\mathcal F}}
\newcommand{\hF}{{\hat{\F}}}
\newcommand{\G}{{\mathcal G}}
\newcommand{\Gp}{{\EuScript Gp}}
\newcommand{\Lie}{{\EuScript Lie}}
\newcommand{\LL}{{\mathcal L}}
\newcommand{\Pa}[1]{$\Pi_{#1}$-algebra}
\newcommand{\PAlg}[1]{\Alg{\Pi_{#1}}}
\newcommand{\Set}{{\EuScript Set}}
\newcommand{\RM}[1]{{#1}\text{-}{\EuScript Mod}}
\newcommand{\SSp}[1]{{\EuScript S}^{#1}_{\F}}
\newcommand{\Ss}{{\mathcal S}}
\newcommand{\Sa}{\Ss_{\ast}}
\newcommand{\sC}[1]{s_{\lra{#1}}\C}
\newcommand{\TT}{{\mathcal T}}
\newcommand{\Ta}{\TT_{\ast}}
%
%
\newcommand{\We}{\mathfrak{W}}
\newcommand{\Cf}{\mathfrak{C}}
\newcommand{\Fb}{\mathfrak{F}}
\newcommand{\Cwcf}{\langle \C; \We,\Cf,\Fb \rangle}
%
%
\newcommand{\bD}{\boldsymbol{\Delta}}
\newcommand{\N}{\mathbb N}
\newcommand{\Q}{\mathbb Q}
\newcommand{\R}{\mathbb R}
\newcommand{\bk}{\mathbf{k}}
\newcommand{\bm}{\mathbf{m}}
\newcommand{\bn}{\mathbf{n}}
\newcommand{\bnp}{\mathbf{n}$+$\mathbf{1}}
\newcommand{\bz}{\mathbf{0}}
\newcommand{\Z}{\mathbb Z}
%
%
\newcommand{\map}{\operatorname{map}}
\newcommand{\B}[2]{\mathbf{B}({#1})_{#2}}
\newcommand{\BGd}[1]{\B{#1}{\bullet}}
\newcommand{\tE}[3]{\tilde{\mathbf{E}}({#1},{#2})_{#3}}
\newcommand{\tEd}[2]{\tE{#1}{#2}{\bullet}}
\newcommand{\EMd}[2]{\mathbf{K}({#1},{#2})_{\bullet}}
\newcommand{\tK}[3]{\tilde{\mathbf{K}}({#1},{#2})_{#3}}
\newcommand{\tKd}[2]{\tK{#1}{#2}{\bullet}}
\newcommand{\M}[1]{\mathbf{M}\lra{#1}}
\newcommand{\Ma}{\M{\alpha}}
\newcommand{\bS}[1]{\mathbf{S}^{#1}}
\newcommand{\X}{\mathbf{X}}
\newcommand{\Y}{\mathbf{Y}}
%
%
\newcommand{\pia}{\pi_{\alpha}}
\newcommand{\pis}{\pi_{\ast}}
\newcommand{\pif}{\pi_{\F}}
\newcommand{\pin}[2]{\hat{\pi}_{#1}{#2}}
%
%
\newcommand{\As}{A_{\ast}}
\newcommand{\Bs}{B_{\ast}}
\newcommand{\Cs}{C_{\ast}}
\newcommand{\Gs}{G_{\ast}}
\newcommand{\Hs}{H_{\ast}}
\newcommand{\Js}{J_{\ast}}
\newcommand{\Ks}{K_{\ast}}
\newcommand{\Ls}{L_{\ast}}
\newcommand{\Ns}{N_{\ast}}
\newcommand{\Ts}{T_{\ast}}
%
%
\newcommand{\Ad}{A_{\bullet}}
\newcommand{\Ads}{A_{\bullet\ast}}
\newcommand{\Bd}{B_{\bullet}}
\newcommand{\Cd}{C_{\bullet}}
\newcommand{\tiEd}{\tilde{E}_{\bullet}}
\newcommand{\Fd}{F_{\bullet}}
\newcommand{\Gd}{G_{\bullet}}
\newcommand{\Kd}{K_{\bullet}}
\newcommand{\tiKd}{\tilde{K}_{\bullet}}
\newcommand{\Pd}{P_{\bullet}}
\newcommand{\Qd}{Q_{\bullet}}
\newcommand{\Qpd}{Q'_{\bullet}}
\newcommand{\Rd}{R_{\bullet}}
\newcommand{\Vd}{V_{\bullet}}
\newcommand{\Wd}{W_{\bullet}}
\newcommand{\Xd}{X_{\bullet}}
\newcommand{\Xdd}{X_{\bullet,\bullet}}
\newcommand{\Yd}{Y_{\bullet}}
\newcommand{\co}[1]{c({#1})_{\bullet}}
\newcommand{\q}[1]{^{({#1})}}
\newcommand{\bd}{\mathbf{d}_{0}}
\newcommand{\dd}[1]{\bd^{#1}}
\newcommand{\bdi}{(\bd)_{\#}}
\newcommand{\ddi}[1]{(\dd{#1})_{\#}}
\newcommand{\jj}[1]{(j_{#1})_{\#}}
\newcommand{\ju}[2]{(j_{#1}^{#2})_{\#}}
%
%
%
\title{Algebraic invariants for homotopy types}
\author{David Blanc}
\address{Dept.\ of Mathematics, University of Haifa, 31905 Haifa, Israel}
\email{blanc@math.haifa.ac.il}
\date{(Revised version) September 18, 1998}
\begin{abstract}
We define a sequence of purely algebraic invariants \ -- \ namely, 
classes in the Quillen cohomology of the \Pa{}\ \ $\pis\X$ \ -- \ for 
distinguishing between different homotopy types of spaces.
Another sequence of such cohomology classes allows one to decide whether a 
given abstract \Pa{}\ can be realized as the homotopy \Pa{}\ of a space.
\end{abstract}

\maketitle
%
%
\sect{Introduction}
\label{ci}

The usual Postnikov system for a (simply-connected) CW complex $\X$ serves to 
determine its homotopy type. One begins with purely algebraic data, 
consisting of the homotopy groups \ $(\pi_{n}\X)_{n=2}^{\infty}$. \ However, 
in order to construct the succesive  approximations \ $\X\q{n}$ \ 
($n\geq 2$), \ with \ $\X\simeq \holim \X\q{n}$, \ one must specify a 
sequence of cohomology classes \ $k_{n}\in H^{n+2}(\X\q{n};\pi_{n+1}\X)$ \ 
(see \cite[IX, \S 2]{GWhE}). These can hardly qualify as \emph{algebraic} 
invariants, since their description involves the cohomology groups of 
topological spaces.
In this paper we show that if one is willing to invest the graded 
group \ $\pis\X\DEF(\pi_{n}\X)_{n=1}^{\infty}$ \ with some further algebraic
structure, the additional information needed to determine the homotopy type of 
$\X$ can be described in purely algebraic terms.

The structure needed on \ $\pis\X$ \ is that of a \Pa{} \ -- \ i.e., a graded 
group equipped with an action of the primary homotopy operations 
(Whitehead products and compositions). \ In this context, the additional
data needed consists of cohomology classes in the Quillen cohomology of this 
\Pa{} \ -- \ which can be defined as usual in algebraic terms 
(see \S \ref{dhom} below). We show:

%
%
\begin{tha}
Given two realizations \ $\X$ \ and \ $\X'$ \ of a \Pa{}\ \ $\Js$, \ there is a 
successively defined sequence of ``difference obstructions'' \ 
$\delta_{n}\in H^{n+1}(\Js,\Omega^{n}\Js)$, \ taking value in the Quillen 
cohomology groups of \ $\Js$, \ with coefficients in the \ $\Js$-module \ 
$\Omega^{n}\Js$, \ whose vanishing implies that \ $\X\simeq\X'$.
\end{tha}

\noindent (See Theorems \ref{ttwo} and \ref{tthree} below).
The \ $(n+1)$-st cohomology class is defined whenever the $n$-th Postnikov 
section of the simplicial space resolutions of the spaces \ $\X$ and $\X'$, \  
respectively, agree, up to homotopy. Even though the obstructions are defined in 
terms of a specific choice of \Pa{}\ resolution of \ $\Js$, \ in fact they depend 
only on the homotopy type of the Postnikov sections\vsm .

Moreover, these cohomology groups can also be used to determine the
realizability of an abstract \Pa{}\ as the homotopy groups of some 
space:

%
%
\begin{thb}
Given a  \Pa{}\ \ $\Js$, \ there is a successively defined sequence of 
``characteristic classes'' \ $\xi\in H^{n+2}(\Js,\Omega^{n}\Js)$, \ which vanish 
if and only if \ $\Js$ \ is realizable by a topological space.
\end{thb}

\noindent (See Theorems \ref{tzero} and \ref{tone} below).
The vanishing requirement should be understood in the sense of 
an obstruction theory: if any such sequence of cohomology classes vanishes,
the \Pa{}\ is realizable; \ if one reaches a non-trivial obstruction, one must 
back-track, and try to vary the choices involved in order to obtain a realization.
These choices again depend only the homotopy type of a suitable Postnikov 
section \ -- \ this time, of a simplicial resolution we are trying to construct 
for the putative topological space $\X$ realizing \ $\Js$. \ 
See Proposition \ref{pfour} below\vsm .

The theory is greatly simplified if we are only interested in the 
\emph{rational} homotopy type of a simply-connected space $\X$. In that case, 
a rational \Pa{}\ is simply a graded Lie algebra over $\Q$, and the cohomology 
theory in question reduces to the usual cohomology of Lie algebras. Theorem A 
thus provides an integral version of (the dual to) the Halperin-Stasheff 
obstruction theory for rational homotopy types  (see \cite{HStaO} and 
\S \ref{rdpr} below).

It is in order to be able to deal with this case, too (and other 
possible variants \ -- \ see \S \ref{ragp} below), that we have stated our 
results for a general model 
catgory $\C$ (subject to certain somewhat restrictive simplifying assumptions 
on $\C$ \ -- \ not all of which are really necessary). 
For technical convenience we have chosen to describe the ordinary topological 
version of our theory within the framework of simplicial groups, rather than 
topological spaces (see \S \ref{srpa} below).

\subsection{notation and conventions}
\label{snac}\stepcounter{thm}

$\TT$ \ will denote the category of topological spaces, and \ $\Ta$ \ that of 
pointed connected topological spaces with base-point preserving maps.
The base-point will be written \ $\ast\in X$.

The category of groups is denoted by \ $\Gp$, \ that of graded groups by \ 
$gr\Gp$, \ that of (left) $R$-modules by \ $\RM{R}$, \ and that of sets by \ 
$\Set$.

\begin{defn}\label{dso}\stepcounter{subsection}
$\bD$ is the category of ordered sequences \ 
$\bn= \langle 0,1,\dotsc,n \rangle$ \ ($n\in \N$), \ with 
order-preserving maps. \ $\bD^{op}$ is the opposite category.
As usual, a \emph{simplicial object} over any category $\C$  is a functor \ 
$X:\bD^{op}\to\C$; \ more explicitly, it is a sequence of objects \ 
$\{ X_{n} \}_{n=0}^{\infty}$ \ in $\C$, equipped with \emph{face maps} \ 
$d_{i}:X_{n}\to X_{n-1}$ \ and \emph{degeneracies} \ 
$s_{j}:X_{n}\to X_{n+1}$ \ ($0\leq i,j\leq n$), \ satisfying the usual
simplicial identities (\cite[\S 1.1]{MayS}). \ We usually denote such a
simplical object by \ $\Xd$. \ The category of simplicial objects over $\C$ is 
denoted by \ $s\C$. \
The standard embedding of categories \ $\co{-}:\C\to s\C$ \ is defined 
by letting \ $\co{X}\in s\C$ \ denote the constant simplicial object \ 
on any \ $X\in\C$ \ (with \ $c(X)_{n}=X$, \ $d_{i}=s_{j}=id_{X}$).
\end{defn}

The category of simplical sets will be denoted by \ $\Ss$, \ rather than \ 
$s\Set$, \ that of pointed connected simplicial sets by \ $\Sa$, \ and that of 
simplicial groups by \ $\G$. \ If we consider a simplicial object \ $\Xd$ \ over 
$\G$, say, we shall sometimes call $n$ in \ $X_{1},\dotsc,X_{n},\dotsc$ \ 
the \emph{external} simplicial dimension, written \ $(-)^{ext}_{n}$, \ 
in distinction from the \emph{internal} simplicial dimension $k$, inside $\G$, 
denoted by \ $(-)^{int}_{k}$. \ In this case we shall sometimes write \ 
$(\Xd)^{int}_{k}\in s\Gp$, \ in contrast with \ $X_{n}\in\G$, \ to emphasize 
the distinction.

The standard $n$ simplex in $\Ss$ is denoted by \ $\Delta[n]$, \ 
generated by \ $\sigma_{n}\in \Delta[n]_{n}$, \ with \ 
$\Lambda^{k}[n]$ \ the subobject generated by \ $d_{i}\sigma_{n}$ \ for \ 
$i\neq k$.

If we denote by \ $\bD\lra{n}$ \ the category obtained from $\bD$ by 
omitting the objects \ $\{\bk\}_{k=n+1}^{\infty}$, \ the category of functors \
$(\bD\lra{n})^{op}\to \C$ \ is called the  category of $n$-\emph{simplicial 
objects} over $\C$ \ -- \ written \ $\sC{n}$. \ If $\C$ has enough 
colimits, the obvious truncation functor \ $\tr{n}:s\C\to \sC{n}$ \ has 
a left adjoint \ $\rho_{n}:\sC{n}\to s\C$, \ and the composite \ 
$\sk{n}\DEF \rho_{n}\circ \tr{n}:s\C\to s\C$ \ is called the 
$n$-\emph{skeleton} functor.

\subsection{organization}
\label{sorg}\stepcounter{thm}

In section \ref{cmc} we review some background material on closed model 
category structures for categories of simplicial objects and show how 
certain convenient CW resolutions may be constructed therein. In section 
\ref{cpf} we construct Postnikov systems for such resolutions, and define
the action of the fundamental group on them; and in section \ref{ccpa} we 
explain how these resolutions are determined in terms of appropriate cohomology 
classes, which may also be used to determine the realizability of a (generalized) 
\Pa{}\ (Theorems \ref{tzero} and \ref{tone}), as well as to distinguish
between different possible realizations (Theorems \ref{ttwo} and \ref{tthree}).

\begin{ack}\label{aa}\stepcounter{subsection}
I would like to thank Dan Kan for suggesting that I continue the project begun 
in \cite{DKStE} and \cite{DKStB}, and Bill Dwyer, Phil Hirschhorn and Emmanuel
Dror-Farjoun for several useful conversations. I am especially grateful to Hans 
Baues and Paul Goerss for pointing out the necessity of taking into account the 
action of the fundamental group in describing the coefficients of the cohomology 
groups. I would also like to thank the referee for his comments.

It should be noted that Baues had previously constructed the first 
difference obstruction of Theorem A, lying in \ $H^{2}(\Js,\Omega\Js)$, \ 
by different methods, and has since extended his construction to the full range 
of invariants we define here: see \cite{BauCF}. Yet a third description of these
invariants, more in the spirit of the original approach of Dwyer, Kan, and Stover,
is planned in \cite{BGoeC}.
\end{ack}
%
%
\sect{model categories of simplicial objects}
\label{cmc}

We first review some background material on model category structures 
for categories of simplicial objects, in particular a slightly expanded 
version of structure defined in \cite{DKStB}, and show how one can construct 
CW resolutions in such a context.

\subsection{model categories}
\label{smc}\stepcounter{thm}

A \emph{model category} in the sense of Quillen (see \cite{QuH}) is a 
category $\C$ equipped with three distinguished classes of morphisms: \ 
$\We$ (weak equivalences), \ $\Cf$, \ and $\Fb$, \ satisfying the following 
assumptions:

\begin{enumerate}
\renewcommand{\labelenumi}{(\arabic{enumi})}
\item \ $\C$ has all small limits and colimits.
\item \ $\We$ \ is a class of quasi-isomorphisms (i.e., there is some 
functor \ $F:\C\to \D$ \ such that \ 
$f\in\We\EQUIV F(f)\text{\ is an isomorphism}$).
\item Any morphism \ $f:A\to B$ \ in $\C$ \ has a factorization \ 
$A\xra{i}C\xra{p}B$ \ ($f=p\circ i$) \ with \ $i\in\Cf\cap \We$ \ and \ 
$p\in\Fb$; \ moreover, this factorization is unique up to weak equivalence, 
in the sense that if \ $A\xra{i'}C'\xra{p'}B$ \ is another such factorization 
of $f$ \ ($i'\in\Cf\cap \We$, \ $p'\in\Fb$), \ then there is a map \ 
$h:C\to C'$ \ such that \ $h\circ i = i'$ \ and \ $p'\circ h=p$.
\item Similarly, any morphism \ $f:A\to B$ \ in $\C$ \ has a factorization \ 
$A\xra{i}C\xra{p}B$ \ ($f=p\circ i$) \ with \ $i\in\Cf$ \ and \ 
$p\in\Fb\cap\We$ \ -- \ again unique up to weak equivalence.
\item We will assume here that the factorizations above may be chosen 
\emph{functorially} (though this is not included in the original definition
in \cite[I, \S 1]{QuH}).
\end{enumerate}

We call the closures under retracts of $\Cf$ and $\Fb$ the classes of 
\textit{cofibrations} and \textit{fibrations}, respectively. \ The definition
given here is then equivalent to the original one of Quillen in \cite{QuH,QuR} \ 
(see \cite[\S 2]{BlaM}).

An object \ $X\in\C$ \ is called \emph{fibrant} if \ $X\to\ast_{f}$ \ is a 
fibration, where \ $\ast_{f}$ \ is the final object of $\C$; similarly $X$ is
\emph{cofibrant} if \ $\ast_{i}\to X$ \ is a cofibration \ 
($\ast_{i}=$ initial object). If \ $X\in\C$ \ 
is cofibrant and \ $Y\in\C$ \ is fibrant, we denote by \ $[X,Y]_{\C}$ \ (or 
simply \ $[X,Y]$) \ the set of homotopy equivalence classes \ $[f]$ \ of maps \
$f:X\to Y$. \ For this to be defined we in fact need only require $X$ to be 
cofibrant \emph{or} $Y$ to be fibrant (cf.\ \cite[I,\S 1]{QuH}).\hsm
A map in \ $\We\cap\Fb$ \ is called a \emph{trivial fibration}, and
one in \ $\We\cap\Cf$ \ a \emph{trivial cofibration}.

Given a model category \ $\langle \C; \We,\Cf,\Fb \rangle$, \ one can 
``invert the weak equivalences'' to obtain the associated \emph{homotopy
category} \ $ho\C$, \ in which the set of morphisms from $X$ to $Y$ is just \ 
$[X,Y]$ \ (at least when $X$ and $Y$ are both fibrant and cofibrant). \ 
See \cite[I]{QuH}, \ \cite[II,\S 1]{QuR}, or \cite[ch.\ IX-XI]{PHirL} for 
some basic properties of model categories.

\subsection{pointed model categories}
\label{spm}\stepcounter{thm}

In a \emph{pointed} model category \ $\langle\C;\We,\Cf,\Fb\rangle$ \ -- \ 
i.e., one with a zero object, denoted by $0$ or $\ast$ 
($=\ast_{f}=\ast_{i}$) \ \ -- \ we may define
the \emph{fiber} of a map (usually: a fibration) \ $f:X\to Y$ \ to be the
pullback of \ $X\xra{f} Y\leftarrow \ast$, \ and the \emph{cofiber} of a 
map (usually: a cofibration) \ $i:A\to B$ \ to be the pushout of \ 
$\ast\leftarrow A\xra{i} B$. \ The \emph{suspension} \ $\Sigma A$ \ of
a (cofibrant) object \ $A\in\C$ \ is then defined to be the cofiber of \ 
$A\amalg A\to A\times I$, \ where \ $A\times I$ \ is any cylinder object for
$A$ \ (cf.\ \cite[I,1,Def.\ 4]{QuH}); \ it is unique up to homotopy 
equivalence. \ Similarly, the \emph{loops} \ $\Omega X$ \ of
a fibrant object $X$ is the fiber of \ $X^{I}\to X\times X$, \ where \ 
$X^{I}$ \ is a path object for $X$ \ (ibid.). 
Finally, the \emph{cone} \ $CA$ \ of a (cofibrant) object \ $A\in\C$ \ is 
the cofiber of either map \ $A\hra A\times I$. \ See \cite[I,2.8-9]{QuH}.

\subsection{simplicial objects}
\label{sso}\stepcounter{thm}

For any category $\C$ \ with coproducts, one has a 
\emph{simplicial structure} (cf.\ \cite[II, \S 1]{QuH}) on the category \ 
$s\C$ \ of simplicial objects over $\C$, defined as usual by:

\begin{enumerate}
\renewcommand{\labelenumi}{(\roman{enumi})}
\item For any simplicial set \ $A\in\Ss$ \ and \ $X\in\C$, \ we define \ 
$X\hotimes A\in s\C$ \ by \ 
$(X\hotimes A)_{n}\DEF \coprod_{a\in A_{n}} X$, \ 
with the face and degeneracy maps induced from those of $A$. \ We 
denote the cofiber of \ $A\hotimes {\ast}\to A\hotimes X$ \ by \ $A\wedge X$.

Now for \ $\Xd\in s\C$ \ we define \ $\Xd\otimes A\in s\C$ \ by \ 
$(\Xd\otimes A)_{n}\DEF \coprod_{a\in A_{n}} X_{n}$ \ (the diagonal of the
bisimplicial object \ $\Xd\hotimes A$).
\item For any \ $\Xd,\Yd\in s\C$ \ we define the function complex \ 
$\map(\Xd,\Yd)$ \ by \ 
$$
\map(\Xd,\Yd)_{n}\DEF\Hom_{s\C}(\Xd\otimes \Delta[n],\Yd),
$$ 
where \ $\Delta[n]\in\Ss$ \ denotes the standard simplicial $n$-simplex.
\end{enumerate}

\begin{defn}\label{dmat}\stepcounter{subsection}
For any complete category $\C$, the \emph{matching object} functor \ 
$M:\Ss^{op}\times s\C\to \C$, \ written \ $M_{A}\Xd$ \ for \ a (finite) 
simplicial set \ $A\in\Ss$ \ and  any \ $\Xd\in s\C$, \ is defined by
requiring that \ $M_{\Delta[n]}\Xd:=X_{n}$, \ and if \ $A=\colim_{i} A_{i}$ \ 
then \ $M_{A}\Xd=\lim_{i} M_{A_{i}}\Xd$ \ (see \cite[\S 2.1]{DKStB}). \ This
may be defined by adjointness, via:
$$
\Hom_{s\C}(Z\otimes A,\Xd)\cong \Hom_{C}(Z,M_{A}\Xd)
$$
\noindent for \ $\Xd\in s\C$ \ and \ $Z\in\C$. 

In particular, we write \ $M^{k}_{n}\Xd$ \ for \ $M_{A}\Xd$ \ where $A$ is the
subcomplex of \ $\sk{n-1}\Delta[n]$ \ generated by the last \ $(n-k+1)$ \ 
faces \ $(d_{k}\sigma_{n},\dotsc,d_{n}\sigma_{n})$. \ When \ $\C=\Set$ \ or \ 
$\Gp$, \ for example, this reduces to:
%
\setcounter{equation}{\value{thm}}\stepcounter{subsection}
\begin{equation}\label{eone}
M^{k}_{n}\Xd=\{(x_{k},\dotsc,x_{n})\in(X_{n-1})^{n+1}~\lvert\ 
d_{i}x_{j}=d_{j-1}x_{i} \text{\ \ for all\ }k\leq i<j\leq n\},
\end{equation}
\setcounter{thm}{\value{equation}}

\noindent and the map \ $\delta^{k}_{n}:X_{n}\to M^{k}_{n}\Xd$ \ induced by 
the inclusion \ $A\hra\Delta[n]$  \ is defined \ 
$\delta_{n}(x)=(d_{k}x,\dotsc,d_{n}x)$. \ The original matching object of
\cite[X,\S 4.5]{BKaH} \ was \ $M_{n}^{0}\Xd=M_{\partial\Delta[n]}\Xd$, \ 
which we shall further abbreviate to \ $M_{n}\Xd$; \ note that each face map \ 
$d_{k}:X_{n+1}\to X_{n}$ \ factors through \ $\delta_{n}\DEF\delta^{0}_{n}$. \ 
See also \S \ref{dpt} below and \cite[XVII, 87.17]{PHirL}.

The dual construction yields the colimit \ $L_{n}\Xd$, \ sometimes called 
the ``$n$-th latching object'' of \ $\Xd$ \ -- \ see 
\cite[\S 2.3(i)]{DKStE}. \ For \ $\Xd\in\S$, \ for example, we have \ 
$L_{n}\Xd\DEF\coprod_{0\leq i\leq n-1} X_{n-1}/\sim$, \ where for any \ 
$x\in X_{n-k-1}$ \ and \ $0\leq i\leq j \leq n-1$ \ we set \ 
$s_{j_{1}}s_{j_{2}}\dotsc s_{j_{k}}x$ \ in the $i$-th copy of \ $X_{n-1}$ \ 
equivalent to \ $s_{i_{1}}s_{i_{2}}\dotsc s_{i_{k}}x$ \ in the $j$-th copy 
of \ $X_{n-1}$ \ whenever the simplicial identity \ 
$$
s_{i}s_{j_{1}}s_{j_{2}}\dotsc s_{j_{k}}=s_{j}s_{i_{1}}s_{i_{2}}\dotsc s_{i_{k}}
$$
\noindent holds (so in particular \ 
$s_{j}x\in (X_{n-1})_{i}$ \ is equivalent to \ $s_{i}x\in (X_{n-1})_{j+1}$ \ 
for all \ $0\leq i\leq j \leq n-1$). \ The map \ 
$\sigma_{n}:L_{n}\Xd\to X_{n}$ \ is defined \ $\sigma_{n}(x)_{i}=s_{i}x$, \ 
where \ $(x)_{i}\in (X_{n-1})_{i}$.
\end{defn}

There are (at least) two ways to extend a given model category structure on 
$\C$ to \ $s\C$:

\begin{defn}\label{dreed}\stepcounter{subsection}
In the \emph{Reedy model structure} on \ $s\C$ \ (see \cite{ReedH} or 
\cite[XVII, \S 88]{PHirL}), \ a simplicial map \ $f:\Xd\to \Yd$ \ is

\begin{enumerate}
\renewcommand{\labelenumi}{(\roman{enumi})}
\item a weak equivalence if \ $f_{n}:X_{n}\to Y_{n}$ \ is a weak equivalence 
in $\C$ \ for each \ $n\geq 0$;
\item a (trivial) cofibration if \ 
$f_{n}\amalg \sigma_{n}:X_{n}\amalg_{L_{n}\Xd}L_{n}\Yd\to Y_{n}$ \ is a 
(trivial) cofibration in $\C$ \ for each \ $n\geq 0$;
\item a (trivial) fibration if \ 
$f_{n}\times \delta_{n}:X_{n}\to Y_{n}\times_{M_{n}\Yd}M_{n}\Xd$ \ is a 
(trivial) fibration in $\C$ \ for each \ $n\geq 0$\vsm .
\end{enumerate}

Note that these definitions imply that \ $\Xd\in s\C$ \ is fibrant if and only
if the maps \ $\delta_{n}:X_{n}\to M_{n}\Xd$ \ are fibrations (in $\C$) for
all $n$.
\end{defn}

We shall require another structure, originally called the ``$E^{2}$-model 
category'' (see \cite[\S 3]{DKStE} and \S \ref{sdrp} below), defined under 
the following

\begin{assume}\label{armc}\stepcounter{subsection}
Assume that \ $\langle \C; \We,\Cf,\Fb \rangle$ \ is a pointed cofibrantly 
generated model category, in which every object is fibrant (this
holds, for example, if \ $\C=\Ta$ \ or \ $\C=\G$). \ Let \ $\F=F_{\C}$ \ be 
a small full subcategory of $\C$ with the following properties:

\begin{enumerate}
\renewcommand{\labelenumi}{(\roman{enumi})}
\item There is a subset \ $\{\Ma\}_{\alpha\in\hF}\subset Obj\F$ \ 
consisting of cogroup objects for $\C$ \ -- \ so there is a natural group 
structure on \ $\Hom_{\C}(\Ma,Y)$ \ for any \ $Y\in\C$. \ 
\item $\F$ is closed under coproducts, and every object \ 
$Z\in\F$ \ is weakly equivalent to some (possibly infinite) coproduct \ 
$\coprod_{i} \M{\alpha_{i}}$ \ with \ $\alpha_{i}\in\hF$ \ -- \ so $Z$ is a 
\emph{homotopy} cogroup object (i.e., \ $[Z,Y]_{\C}$ \ has a natural group 
structure). However, we do not require the morphisms in $\F$ to respect the 
cogroup structure, even up to homotopy.
\item  $\F$ is closed under suspensions \ -- \ that is, for each \ 
$X\in\F$, \ there is a model for \ $\Sigma X$ \ in $\F$. \ 
We also assume \ $C\Ma\in \F$ \ for every \ $\alpha\in\hF$ \ 
(\S \ref{spm}).
\end{enumerate}
\end{assume}

We now wish to define an algebraic model for the collection of sets of 
homotopy classes of maps \ $\{[X,Y]_{\C}\}_{X\in\F}$, \ for a given object \ 
$Y\in\C$. \ This is provided by the following

\begin{defn}\label{dfpa}\stepcounter{subsection}
Given \ $\F\subset\C$ \ as in \S \ref{armc}, we define a \ 
\emph{\Pa{\F}}\ to be a functor \ $ho(\F)^{op}\to\Set$, \ which takes 
coproducts in $\F$ to products in \ $\Set$ \ (compare \cite{DrecP}).

The category of all \Pa{\F}s will be denoted by \ $\PAlg{\F}$, \ and 
the functor \ $[ho\F,-]:\C\to \PAlg{\F}$ \ defined \ 
$([\BB,Y])_{\BB\in ho\F}$ \ will be denoted by \ $\pif$. \ $\PAlg{\F}$ \ is a
category of universal graded algebras, or CUGA, in the sense of 
\cite[\S 2.1]{BStG}. \ In particular, the \emph{free} \Pa{\F}s are
those isomorphic to \ $\pif X$ \ for some \ $X\in \F$. \ If we assume that \ 
$X\simeq \coprod_{\alpha\in\hF}\coprod_{t\in T_{\alpha}}\Ma_{t}$ \ for some \ 
$\hF$-graded set \ $\Ts$, \ we say that \ $\pif X$ \ is the free \Pa{\F} 
\emph{generated by} \ $\Ts$.

If \ $f:X\to Y$ \ is a morphism in \ $\C$, \ the induced morphism of
\Pa{\F}s, \ $\pif f:\pif X\to\pif Y$, \ will be denoted simply by \ $f_{\#}$.

\end{defn}

\begin{remark}\label{rfpa}\stepcounter{subsection}
Since all objects in $\F$ are homotopy equivalent to coproducts of objects
from the set \ $\hF$, \ a \Pa{\F}\ may be thought of 
more concretely as an \ $\hF$-graded group \ -- \ i.e., a collection  of 
groups \ $(G_{\alpha})_{\alpha\in\hF}$ \ -- \ equipped with a 
(contravariant) action of the homotopy classes of morphisms in $\F$ on them,
modeled on the action of such homotopy classes on \ 
$\{[\Ma,Y]\}_{\alpha\in\hF}$ \ by precomposition (cf.\ \cite[XI, \S 1]{GWhE}).

We shall write \ $\pia X$ \ for \ $(\pif X)_{\alpha}\DEF[\Ma,X]$, \ 
and \ $\pi_{\alpha+k}X$ \ for \ $[\Sigma^{k}\Ma,X]$.
\end{remark}

\begin{defn}\label{dapa}\stepcounter{subsection}
As usual, a \Pa{\F}\ $X$ is called \emph{abelian} if \ 
$\Hom_{\PAlg{\F}}(X,A)$ \ has a natural abelian group structure for any \ 
$A\in\PAlg{\F}$ \ (see \cite[\S 5.1]{BStG} for an explicit description.).
In particular, for any \ $X\in\PAlg{\F}$, \ its \emph{abelianization} \ 
$X_{ab}$ \ may be defined as in \cite[\S 5.1.4]{BStG} as a suitable quotient 
of $X$. Another abelian \Pa{\F}\ which may be defined for any $X$ is its
\emph{loop algebra} \ $\Omega X$, \ defined by \ 
$\Omega X(\BB)\DEF X(\Sigma\BB)$ \ (cf.\ \cite[\S 9.4]{DKStB}; recall that 
$\F$ is closed under suspension). The fact that it is abelian follows as in 
\cite[Prop.\ 9.9]{GrayH}. \ 
The (abelian) category of abelian \Pa{\F}s will be denoted by \ $\PAlg{\F}_{ab}$.
\end{defn}

\begin{example}\label{erep}\stepcounter{subsection}
In \ $\C=\Ta$, \ let \ $\F$ \ denote the subcategory whose objects are wedges
of spheres of various dimensions; then for any space \ $X\in\Ta$, \ 
the functor \ $\pif X$ \ is determined up to isomorphism by \ $\pis X$, \ 
the homotopy \emph{\Pa{}} of $X$ \ -- \ that is, its homotopy
groups, together with the action of the primary homotopy operations (Whitehead
products and compositions) on them. See \cite[\S 2]{BlaA} or \cite[\S 4]{StoV}.
In particular, the abelian \Pa{}s are those for which all Whitehead products
are trivial (cf.\ \cite[\S 3]{BlaA}).
\end{example}

\begin{remark}\label{rgp}\stepcounter{subsection}
This example does not quite fit our assumptions (\S \ref{armc}), \ since the
spheres are only co-$H$-spaces, i.e., \emph{homotopy} cogroup objects 
in \ $\Ta$. \ This does not affect the arguments at this stage \ -- \ in fact,
this is the original example of an \ ``$E^{2}$-model category'' in 
\cite{DKStE}. However, for our purposes $\G$ appears to be more convenient 
than \ $\Ta$ \ as a model for the homotopy category of (connected) spaces  \ 
(see \cite{KanT}; also, e.g., \cite[\S 5]{BlaL}). 

In fact, in all the examples we have in mind the objects in $\C$ will have an 
(underlying) group structure, so it will be convenient to add to \S \ref{armc} 
the following additional
\end{remark}

\begin{assume}\label{agp}\stepcounter{subsection}
$\C$ is equipped with a faithful forgetful functor \ $\hat{U}:\C\to\D$ \ -- \ 
where \ $\D$ is one of the ``categories of groups'' \ $\D=\Gp$, \ $gr\Gp$, \ 
$\G$, \ $\RM{R}$, \ or \ $s\RM{R}$, \ for some ring $R$ \ -- \ and the cogroup 
objects \ $\Ma\in\hF$ \ of \S \ref{armc}(i) are in the image of its adjoint 
$\hat{F}$, \ with the group structure on \ $\Hom_{\C}(\Ma,X)$ \ induced from 
that of \ $\hat{U}(X)$. \ When \ $\D=\G$ \ or \ $\D=s\RM{R}$, \  the objects \ 
$\Ma$ \ must actually lie in the image of the composite \ 
$\hat{F}\circ F':\Ss\to \C$, \ where \ $F':\Ss\to\D$ \ is adjoint to the 
forgetful functor \ $U':\D\to\Ss$.

We also assume that the adjoint pair \ $(\hat{U},\hat{F})$ \ \emph{create} the 
model category structure on $\C$ in the sense of \cite[\S 4.13]{BlaM} \ -- \  
so in particular $\hat{U}$ creates all limits in $\C$ 
(cf.\ \cite[V,\S 1]{MacC}).
\end{assume}

\begin{remark}\label{ragp}\stepcounter{subsection}
In fact, the categories $\C$ in which shall be interested are the following:
\begin{enumerate}
\renewcommand{\labelenumi}{$\bullet$}
\item $\C=\G$, \ so \ $s\C$, \ the category of bisimplicial groups, is a 
model for simplicial spaces; 
\item $\C=\Gp$, \ so \ $s\C=\G$ \ is a model for the homotopy category of
connected topological spaces of the homotopy type of a CW complex;
\item $\C=d\LL$, \ the category of differential graded Lie algebras (or 
equivalently, \ $\C=s\Lie$), \ so \ $s\C$ \ is a model for simplicial rational 
spaces;
\item $\C=\Lie$, \ the category of Lie algebras, so \ $s\Lie$ \ is a model
for (simply connected) rational spaces \ (cf.\ \cite[II,\S 4-5]{QuR});
\item $\C=\RM{R}$, \ the category of (left) modules over a not-necessarily
commutative, possibly graded, ring $R$, \ so \ $s\C$ \ is a model for chain 
complexes over $R$.
\end{enumerate}

\noindent and it is the desire to give a unified treatment for these
five cases that forces upon us the somewhat unnatural set of assumptions we
have made in \S \ref{armc} and here.
\end{remark}

\begin{defn}\label{dffm}\stepcounter{subsection}
A map \ $f:\Vd\to\Yd$ \ in  \ $s\C$ \ is called \emph{$\F$-free} if for 
each \ $n\geq 0$, \ there is 

\begin{enumerate}
\renewcommand{\labelenumi}{\alph{enumi})\ }
\item a cofibrant object \ $W_{n}$ \ which is weakly equivalent to an 
object in $\F$;
\item a map \ $\varphi_{n}:W_{n}\to Y_{n}$ \ in $\C$ which induces a trivial 
cofibration \ $(V_{n}\amalg_{L_{n}\Vd}L_{n}\Yd)\amalg W_{n}\to Y_{n}$.
\end{enumerate}
\end{defn}

\subsection{The resolution model category}
\label{srmc}\stepcounter{thm}

Given a model category $\C$ and a subcategory $\F$ as in \S \ref{armc},
we define the \emph{resolution model category structure} on \ $s\C$, \ 
with respect to $\F$ by setting a simplicial map \ $f:\Xd\to \Yd$ \ to be

\begin{enumerate}
\renewcommand{\labelenumi}{(\roman{enumi})}
\item a \emph{weak equivalence} if \ $\pif f$ \ is a weak equivalence
of $\hF$-graded simplicial groups (\S \ref{rfpa}).
\item a \emph{cofibration} if it is a retract of an $\F$-free map;
\item a \emph{fibration} if it is a Reedy fibration (Def.\ \ref{dreed}(iii)) 
and \ $\pif f$ \ is a (levelwise) fibration of simplicial groups 
(that is, for each \ $B\in\F$ \ and each \ $n\geq 0$, \ the group 
homomorphism \ $[B,X_{n}]\xra{[B,f_{n}]}[B,Y_{n}]^{ext}$ \ is an epimorphism \ 
(where for  \ $\Gd:=[B,\Yd]\in\G$, \ $\Gd^{ext}$ \ denotes the connected 
component of the identity) \ -- \ see \cite[II,3.8]{QuH}.
\end{enumerate}

This was originally called the ``$E^{2}$-model category structure'' on \ 
$s\C$. \ See \cite[\S 5]{DKStE} for further details. 

\begin{example}\label{ermc}\stepcounter{subsection}
Let \ $\C=\Gp$ \ with the \emph{trivial} model category structure: \ i.e.,
only isomorphisms are weak equivalences, and every map is both a fibration and
a cofibration. Let \ $\F_{\Gp}$ \ be the category of all free groups 
(which are the cogroup objects in \ $\Gp$ \ -- \ cf.\ \cite{KanMD}). \ The 
resulting resolution model category structure on \ $\G\DEF s\Gp$ \ is the usual 
one (cf.\ \cite[II, \S 3]{QuH}). \ This observation is due to Pete Bousfield.
We can then iterate the process by letting \ 
$\F_{\G}$ \ be the category of (coproducts of) the $\G$-spheres, defined: \ 
$\bS{n}\DEF FS^{n-1}\in\G$ \ -- \ see \cite{MilnF} \ -- \ (with \ 
$\bS{0}\DEF GS^{0}$), and obtain a resolution model category structure on \ 
$s\G$ \ (bisimplicial groups). 

Note that if we tried to do the same for \ $\C=\Set$, \ there are no 
nontrivial cogroup objects, while in $\Ss$ not all objects are fibrant 
(see \S \ref{armc}). \ The category \ $\Ta$ \ of pointed topological spaces, 
which is the main example we actually have in mind, does not quite fit our 
assumptions (but see \S \ref{rgp} above).
\end{example}

Motivation for the name of ``resolution model category'' is provided by 
the following

\begin{defn}\label{dres}\stepcounter{subsection}
A \emph{resolution} of an object \ $\Xd\in s\C$ \ (relative to $\F$) 
is a \emph{cofibrant replacement} for \ $\Xd$ \ in the resolution 
model category on \ $s\C$ \ determined by $\F$: that is, it is any
cofibrant object \ $\Qd$, \ equipped with a weak equivalence to \ $\Xd$, \ 
which may be obtained from the factorization of \ $\ast\to\Xd$ \ as \ 
$\ast\xra{\text{cof}}\Qd\xra{\text{fib+w.e.}}\Xd$ \ -- \ 
and is thus unique up to weak equivalence, by \S \ref{smc}(4)).

More classically, a (simplicial) resolution for an object \ $X\in\C$ \ is 
a resolution of the constant simplicial object \ $\co{X}$ \ 
(cf.\ \S \ref{dso}) in \ $s\C$.
\end{defn}

\subsection{Functorial resolutions}
\label{sfres}\stepcounter{thm}

The construction of \cite[\S 2]{StoV} provides canonical resolutions in \ 
$s\C$, \ defined as follows: consider the comonad \ $L:\C\to\C$ \ given by
\setcounter{equation}{\value{thm}}\stepcounter{subsection} 
\begin{equation}\label{etwo}
 LY = \ \coprod _{\alpha\in\hF} \ 
         \coprod _{\phi\in\Hom_{\C}(\Ma,Y)} \ \Ma_{\phi} \ \ 
         \bigcup \ \ \coprod _{\alpha\in\hF} \ 
         \coprod _{\Phi\in\Hom_{\C} (C\Ma,Y)} \ C\Ma_{\Phi},
\end{equation}
\setcounter{thm}{\value{equation}}

\noindent by which we mean the the coproduct, over all \ 
$\phi:\Ma\to Y$, \ of the colimits of the various diagrams consisting 
of an inclusion \ $\Ma_{\phi}\to C\Ma_{\Phi}$ \ for each \ 
$\Phi:C\Ma\to Y$ \ such that \ $\Phi\rest{\Ma}=\phi$. \ 
The counit \ $\varepsilon:LY\to Y$ \ is ``evaluation of indices'', and the 
comultiplication \ $\vartheta:LY\hra L^{2}Y$ \ is the obvious ``tautological'' 
one. Note that \ $LY\in\F$ \ for any \ $Y\in\C$ \ by our assumptions on $\F$ 
(\S \ref{armc}).

Given \ $X\in\C$, \ we define its \emph{canonical resolution} \ 
$\Qd\to X$ \ by \ $Q_{n}\DEF L^{n+1}X$, \ with the degeneracies and face maps 
induced as usual by \ $\varepsilon$ and $\vartheta$ \ 
(see \cite[App., \S 3]{GodT}).

The construction can be modified so as to yield
resolutions for arbitrary \ $\Yd\in s\C$, \ and not only \ $\co{X}$. \ 
Moreover, it has the advantage that \ $\pif g:\pif Q_{n}\to \pif Y_{n}$ \ is
clearly surjective for all $n$, \ so $g$ can be changed into a fibration
(Def.\ \ref{srmc}(iii)) by simply changing each \ $Q_{n}$ \ up to homotopy,
which yields the factorization needed for \S \ref{smc}(3).

An alternative (noncanonical) construction of a resolution is given in 
Proposition \ref{ptwo} below.

\subsection{representing objects for $s\C$}
\label{sros}\stepcounter{thm}

Just as the spheres ``represent'' the weak equivalences in the usual model
structure on \ $\Ta$, \ for example, in the sense that a map \ $f:X\to Y$ \ 
is a weak equivalence if and only if it induces an isomorphism \ 
$f_{\ast}:[S^{n},X]\to[S^{n},Y]$ \ for each \ $n\geq 0$, \ we may similarly 
define representing objects for the resolution model category (compare 
\cite[\S 5.1]{DKStB}):

\begin{defn}\label{dsfs}\stepcounter{subsection}
Given a model category $\C$ and a subcategory $\F$ as above,
for each \ $n\geq 0$, \ the \emph{$n$-dimensional simplicial $\F$-sphere}, \ 
denoted by \ $\SSp{n}$, \ is the subcategory \ $\Sigma^{n}\F$ \ of \ 
$s\C$, \ whose objects are of the form \ 
$\Sigma^{n}X\DEF X\wedge S^{n}$ \ for \ $X\in\F$, \ where \ 
$S^{n}=\Delta[n]/\dot{\Delta}[n]$ \ is the usual
simplicial $n$-sphere (see \S \ref{sso}(i)).

Note that each such \ $\Sigma^{n}X$ \ is cofibrant (in fact, free) in the 
resolution model category \ $s\C$. \ Moreover, by the definition
of the simplicial structure on \ $s\C$ \ (\S \ref{sso}), \ $\Sigma^{n}X$ \ 
is also a cogroup object in \ $s\C$.

Given \ $\Yd\in s\C$, \ choose some \emph{fibrant replacement} \ $\Xd$ \ 
(that is, factor \ $\Yd\to\ast$ \ as \ 
$\Yd\xra{\text{cof+w.e.}}\Xd\xra{\text{fib}}\ast$, \ using \S \ref{smc}(3)) \ 
and define \ $\pin{n}{\Yd}$ \ (also written \ $[\SSp{n},\Yd]$) \ to be the 
$\hF$-graded set \ $\pi_{0}\map(\SSp{n},\Xd)$. \ 
This definition is independent of the choice of $\Xd$.

We define a map \ $f:\Xd\to\Yd$ \ in \ $s\C$ \ to be an 
\emph{$\F$-equivalence} \ if it induces isomorphisms in \ $\pin{n}{(-)}$ \ 
for all \ $n\geq 0$. \ 
\end{defn}

\subsection{fibration sequences}
\label{sfs}\stepcounter{thm}

Let \ $\F\subset\C$ \ be as in \S \ref{armc}, \ and \ $\Xd\to\Yd$ \ a 
fibration in the resolution model category \ $s\C$ \ (\S \ref{srmc}), with
fiber \ $\Fd$ \ (\S \ref{spm}). \ Then as usual we have the 
\emph{long exact sequence of the fibration}:
%
\setcounter{equation}{\value{thm}}\stepcounter{subsection}
\begin{equation}\label{ethree}
\dotsb \to \pin{n+1}{\Yd} \xra{\partial_{\ast}} \pin{n}{\Fd} \to 
\pin{n}{\Xd} \to \pin{n}{\Yd} \to \dotsb \to \pin{0}{\Yd},
\end{equation}
\setcounter{thm}{\value{equation}}

\noindent (see \cite[I,3.8]{QuH}), which in fact may be constructed in 
this case as for \ $\Sa$ \ (see \cite[7.6]{MayS}).

\begin{defn}\label{dnc}\stepcounter{subsection}
Given \ $\Xd\in s\C$, \ , we define the 
$n$-\emph{cycles} object of \ $\Xd$, \ written \ $Z_{n}\Xd$, \ to be the 
fiber of \ $\delta_{n}:X_{n}\to M_{n}\Xd$ \ (see \S \ref{dmat}), so \ 
$Z_{n}\Xd=\{ x\in X_{n}\,|\ d_{i}x=0 \ \text{for}\ i=0,\dotsc,n\}$ \ 
(cf.\ \cite[I,\S 2]{QuH}). \ 
Of course, this definition really makes sense only when \ $\delta_{n}$ \ is
a fibration in $\C$. \ Similarly, the $n$-\emph{chains} object
of \ $\Xd$, \ written \ $C_{n}\Xd$, \ is defined to be the fiber of \ 
$\delta^{1}_{n}:X_{n}\to M^{1}_{n}\Xd$.
\end{defn}

Note that for any \ $W\in\C$ \ and fibrant \ $\Xd\in s\C$ \ we have natural 
adjunction isomorphisms \ 
$\Hom_{s\C}(W\wedge S^{n},\Xd)\cong \Hom_{\C}(W,Z_{n}\Xd)$ \ and \ 
$\Hom_{s\C}(W\wedge D^{n},\Xd)\cong \Hom_{\C}(W,C_{n}\Xd)$ \ (where \ 
$D^{n}\DEF \Delta^{n}/\Lambda^{0}[n]\in \Ss$ \ is a simplicial model for the
$n$-disc)\vsm .

If \ $\Xd$ \ is fibrant, the map \ 
$\bd=\dd{n}\DEF d_{0}\rest{C_{n}\Xd}:C_{n}\Xd\to Z_{n-1}\Xd$ \ is the 
pullback of \ $\delta_{n}:X_{n}\to M_{n}\Xd$ \ along the inclusion \ 
$\iota:Z_{n-1}\Xd\to M_{n}\Xd$ \ (where \ $\iota(z)=(z,0,\dotsc,0)$), \ so \ 
$\bd$ \ is a fibration (in $\C$), \ fitting into a fibration sequence
%
\setcounter{equation}{\value{thm}}\stepcounter{subsection}
\begin{equation}\label{efour}
\dotsb \Omega Z_{n-1}\Xd\to Z_{n}\Xd\xra{j^{X}_{n}}C_{n}\Xd \xra{\bd} 
Z_{n-1}\Xd
\end{equation}
\setcounter{thm}{\value{equation}}

\noindent (see \cite[Prop.\ 5.7]{DKStB}).  \ Moreover, there is an exact 
sequence of \Pa{\F}s \ 
%
\setcounter{equation}{\value{thm}}\stepcounter{subsection}
\begin{equation}\label{efive}
\pif C_{n+1}\Xd \xra{\bdi} \pif Z_{n}\Xd \xra{q} \pin{n}{\Xd}\to 0,
\end{equation}
\setcounter{thm}{\value{equation}}

\noindent (see \cite[Prop.\ 5.8]{DKStB}), \ which provides a (relatively) 
explicit way to recover \ $\pin{n}{\Xd}$ \ from \ $\Xd$.

Finally, the composition of the boundary map \ 
$\partial_{\ast}:\Omega Z_{n-1}\Xd \to Z_{n}\Xd$ \ of the fibration sequence 
\eqref{efour} with \ $\Omega\bd$ \ is trivial, so by \eqref{efive} it induces 
a map of \Pa{\F}s from \ $\pin{n-1}{\Omega\Xd} \cong \Omega\pin{n-1}{\Xd}$ \ 
(\S \ref{dapa}) to \ $\pif Z_{n}\Xd$ \ which, composed with the  
map $q$ in \eqref{efive}, defines a ``shift map'' \ 
$s:\Omega\pin{n-1}{\Xd}\to \pin{n}{\Xd}$ \ (see \cite[Prop.\ 6.2]{DKStB}).

\subsection{the simplicial \Pa{\F}}
\label{ssp}\stepcounter{thm}

Applying the functor \ $\pif$ \ dimensionwise to any simplicial object \ 
$\Xd\in s\C$ \ yields a simplicial \Pa{\F} \ $\Gd=\pif\Xd$, \ which is in 
particular an $\hF$-graded simplicial group; its homotopy groups form a 
sequence of $\hF$-graded groups which we denote by \ 
$(\pi_{n}\pif\Xd)_{n=0}^{\infty}$, \ and each \ $\pi_{n}\pif\Xd$ \ is a 
\Pa{\F}.

Note that as for any (graded) simplicial group, the homotopy groups of \ 
$\Gd$ \ may be computed using the Moore chains \ $\Cs\Gd$, \ defined \ 
$C_{n}\Gd\DEF \cap_{i=0}^{n}\Ker\{d_{i}:G_{n}\to G_{n-1}\}$ \ (cf.\ 
\S \ref{dnc} and \cite[17.3]{MayS}), \ and we have the following version of
\cite[Prop.\ 2.11]{BlaCW}

%
%
\begin{lemma}\label{lone}\stepcounter{subsection}
For any fibrant \ $\Xd\in s\C$, \ the inclusion \ $\iota:C_{n}\Xd\hra X_{n}$ \ 
induces an isomorphism \ 
$\iota_{\star}:\pis C_{n}\Xd\cong C_{n}(\pis\Xd)$ \ for each \ $n\geq 0$.
\end{lemma}

\begin{proof}
(a)\hs First, note that any \emph{trivial} cofibration \ $j:A\hra B$ \ 
in $\Ss$  \ induces a fibration \ $j^{\ast}:M_{B}\Xd\to M_{A}\Xd$ \ in $\C$.

To see this, by assumption \ref{agp} it suffices to consider \ $\C=\D$, \ 
(since by \cite[Def.\ 4.13]{BlaM}, \ $f$ is a fibration in $\C$ if and 
only if \ $Uf$ \ is a fibration in $\D$), and in fact the only nontrivial 
case is when \ $\D=\G$ \ (where the fibrations are maps which surject onto 
the identity component \ -- \ see \cite[II, 3.8]{QuH}). \ Note that in internal
simplicial dimension $k$ we have \ 
$(M_{A}\Xd)^{int}_{k}\cong\Hom_{s\Gp}(FA,(\Xd)^{int}_{k})$ \ 
(see \S \ref{dso}) for \ $A\in\Ss$, \ where \ $F$ \ denotes the 
(dimensionwise) free group functor. 
Since \ $FA$ \ is fibrant in \ $s\Gp$, \ $Fj:FA\hra FB$ \ has a left inverse \ 
$r:FB\to FA$, \ so \ $j^{\ast}:(M_{B}\Xd)^{int}_{k}\to (M_{A}\Xd)^{int}_{k}$ \ 
has a right inverse \ $r^{\ast}$, \ so in particular is onto.  Since this
is true in each simplicial dimension $k$, \ $j^{\ast}:M_{B}\Xd\to M_{A}\Xd$ \ 
is a fibration in $\G$. (Note that \ $d_{i}:X_{n}\to X_{n-1}$ \ is always a 
fibration\vsm .)

(b)\hs In addition, \ $\psi^{k}_{n}=j^{\ast}:M^{0}_{n}\Xd\to M^{k}_{n}\Xd$ \ 
is a fibration for all \ $0\leq k\leq n$, as one can see by considering 
\eqref{eone} (since \ $\delta_{n-1}$ \ surjects onto the identity 
component by assumption)\vsm . 

(c)\hs Given \ $\eta\in C_{n}(\pia \Xd)$, \ represented by 
$h:\Ma\to Y_{n}$, \ with \ $d_{j}h\sim 0$ \ for \ $1\leq j\leq n$, \ note 
that for \ $1\leq k\leq n$, \  $M^{k}_{n}\Xd$ \ is the pullback of \ 
$$
M^{k+1}_{n}\Xd\xra{(d_{k},\dotsc,d_{k})} M^{k}_{n-1}\Xd \xla{\delta^{k}_{n-1}}
X_{n-1},
$$

\noindent in which \ $(d_{k},\dotsc,d_{k})$ \ is a fibration by (a) if \ 
$k\geq 1$, \ so this is in fact a homotopy pullback square 
(see \cite[\S 1]{MatP}).  \ By descending induction on \ $1\leq k\leq n$, \ 
(starting with \ $\delta^{n}_{n}=d_{n}$), \ we may assume \ 
$\delta^{k+1}_{n}\circ h:\Ma\to M^{k+1}_{n}\Xd$ \ is nullhomotopic in $\C$,
as is \ $d_{k}\circ h$, \ so the induced pullback map, which is just \ 
$\delta^{k}_{n}\circ h:\Ma\to M^{k}_{n}\Xd$, \ is also nullhomotopic by the
universal property. We conclude that \ $\delta^{1}_{n}\circ h\sim 0$, \ 
and since \ $\delta^{1}_{n}:X_{n}\to M^{1}_{n}\Xd$ \ is a fibration by (b), 
we can replace $h$ by a homotopic map \ $h':\Ma\to X_{n}$ \ such that \ 
$\delta_{n} h'=0$. \ Thus \ $h'$ \ lifts to \ $Z_{n}\Yd=\Fib(\delta_{n})$, \ 
so \ $\iota_{\star}$ \ is surjective\vsm .

(d)\hs Even though the retraction \ $r:F\Delta[n]\to F\Lambda^{0}_{n}$ \ in (a)
is not canonical, it may be chosen independently of the internal simplicial
dimension $k$ to yield a section \ $r^{\ast}$ \ for \ 
$\delta^{1}_{n}=j^{\ast}:X_{n}\epic M^{1}_{n}\Xd$. \ 
The long exact sequence in \ $[\Ma,-]$ \ for the fibration 
sequence \ $C_{n}\Yd\xra{i} Y_{n}\xra{\delta'_{n}} M^{1}_{n}\Yd$ \ 
(cf.\ \cite[I,\S 3]{QuH}) then implies that \ $i_{\#}$ \ is monic, so \ 
$\iota_{\star}$ \ is, too. The argument lifts from \ $\D=\G$ \ to $\C$ because 
the objects \ $\Ma$ \ are in the image of the adjoint of \ $U:\C\to\D$, \ by
assumption \ref{agp}.
\end{proof}

This Lemma, together with \eqref{efive}, yields a commuting diagram:

%
%
\begin{figure}[htbp]
\begin{picture}(100,90)(-110,-10)
\put(5,60){$\pif C_{n+1}\Xd$}
\put(62,65){\vector(1,0){70}}
\put(85,70){$\bdi$}
\put(135,60){$\pif Z_{n}\Xd$}
\put(180,65){\vector(1,0){38}}
\put(218,65){\vector(1,0){4}}
\put(225,60){$\pin{n}{\Xd}$}
\put(30,55){\vector(0,-1){43}}
\put(17,32){$\iota_{\star}$}
\put(33,32){$\cong$}
\put(0,0){$C_{n+1}(\pif\Xd)$}
\put(65,5){\vector(1,0){57}}
\put(80,10){$d_{0}^{\pif\Xd}$}
\put(150,55){\vector(0,-1){43}}
\put(137,32){$\hat{\iota}_{\star}$}
\put(125,0){$Z_{n}(\pif \Xd)$}
\put(180,5){\vector(1,0){33}}
\put(213,5){\vector(1,0){4}}
\put(220,0){$\pi_{n}\pif\Xd$}
\multiput(235,55)(0,-3){14}{\circle*{.5}}
\put(235,12){\vector(0,-1){3}}
\put(240,32){$h$}
\end{picture}
\caption[fig1]{}
\label{fig1}
\end{figure}

\noindent which defines the dotted morphism of \Pa{\F}s \ 
$h:\pin{n}{\Xd}\to\pi_{n}(\pif\Xd)$ \ (this was called the ``Hurewicz map'' 
in \cite[7.1]{DKStB}). Note that for \ $n=0$ \ the map \ 
$\hat{\iota}_{\star}$ \ is an isomorphism, so $h$ is, too.

\subsection{An exact couple}
\label{sec}\stepcounter{thm}

If \ $\Xd\in s\C$ \ is Reedy fibrant, the long exact sequences 
\eqref{ethree} for the fibrations \ $\C_{n+1}\Xd\to Z_{n}\Xd$ \ fit into an \ 
$(\N,\hF)$-bigraded exact couple \ 
$(D^{1}_{\ast,\alpha},E^{1}_{\ast,\alpha})$ \ with \ 
$D^{1}_{k,\alpha}\cong \pia Z_{k}\Xd$ \ and \ 
$E^{1}_{k,\alpha}\cong \pia C_{k}\Xd$ \ for \ $k\geq 0$ \ and \ 
$\Ma\in\hat{F}$. \ As in \cite[\S 8]{DKStB} \ the derived couple has \ 
$D^{2}_{k,\alpha}\cong (\pin{k}{\Xd})_{\alpha}$ \ and \ 
$E^{2}_{k,\alpha}\cong \pi_{k}(\pia \Xd)$ \ (using Lemma \ref{lone}), 
which fit into a ``spiral exact sequence'' 
%
\setcounter{equation}{\value{thm}}\stepcounter{subsection}
\begin{equation}\label{esix}
\dotsb \to \pi_{n+1}\pif\Xd \xra{\partial} \Omega\pin{n-1}{\Xd} \xra{s} 
\pin{n}{\Xd} \xra{h} \pi_{n}\pif\Xd \to \dotsb 
\pin{0}{\Xd}\xra{h}\pi_{0}\pif\Xd \to 0
\end{equation}
\setcounter{thm}{\value{equation}}

\noindent as in \cite[8.1]{DKStB}, so by Reedy fibrant replacement 
(\S \ref{dsfs}), one has such an exact sequence for any \ $\Yd\in s\C$. \ 
Of course, \ $\pin{-1}{\Xd}\DEF 0$; \ and at the right hand end we have \ 
$h:\pin{0}{\Xd}\cong\pi_{0}\pif\Xd$, \ as noted above.

We immediately deduce the following
%
%
\begin{prop}\label{pone}\stepcounter{subsection}
A map \ $f:\Xd\to\Yd$ \ in \ $s\C$ \ is a weak equivalence in the 
resolution model category \ -- \ i.e., induces an isorphism in \ 
$\pi_{n}\pif$ \ for all \ $n\geq 0$ \ (\S \ref{srmc}(i)) \ -- \ 
 if and only if it is an $\F$-equivalence \ -- \ i.e., \ induces an 
isomorphism in \ $\pin{n}{}$ \ for all \ $n\geq 0$ \ (see \S \ref{dsfs}).
\end{prop}

\subsection{Resolutions}
\label{sres}\stepcounter{thm}
By Definition \ref{dres}, a \emph{resolution} of an object \ $X\in\C$ \ is
a simplicial object \ $\Qd$ \ over $\C$ which is cofibrant and has a weak 
equivalence \ $f:\Qd\to \co{X}$. \ Note that such an $f$ is detemined by 
an \emph{augmentation} \ $\varepsilon:Q_{0}\to X$ \ in \ $\C$ \ (with \ 
$d_{0}\circ \varepsilon=d_{1}\circ \varepsilon$); \ by Proposition \ref{pone}, 
$f$ is a weak equivalence if and only if the augmented 
$\hF$-graded simplicial group \ $\varepsilon_{\ast}:\pif\Qd\to \pif X$ \ 
is acyclic (i.e., has vanishing homotopy groups in all dimensions $\geq 0$).

The long exact sequence \eqref{esix} then implies that \ 
%
\setcounter{equation}{\value{thm}}\stepcounter{subsection}
\begin{equation}\label{eseven}
\pin{n}{\Qd}\cong \Omega^{n}\pif X \ \ \ \ \text{for all }n\geq 0.
\end{equation}
\setcounter{thm}{\value{equation}}

\begin{defn}\label{dcw}\stepcounter{subsection}
A \emph{CW complex} over a pointed category $\C$ is a simplicial object \ 
$\Rd\in s\C$, \ together with a sequence of objects \ $\bar{R}_{n}$ \ 
($n=0,1,\dotsc$) \ -- called a \emph{CW basis} for \ $\Rd$ \ -- \ such that \ 
$R_{n}=\bar{R}_{n}\amalg L_{n}\Rd$ \ (\S \ref{dmat}), \ and \ 
$d_{i}\rest{\bar{R}_{n}}=0$ \ for \ $1\leq i\leq n$. \ The morphism \ 
$\bar{d}^{n}_{0}:\bar{R}_{n}\to Z_{n-1}\Rd$ \ is called the $n$-th
\emph{attaching map} for \ $\Rd$ \ (compare \cite[\S 5]{BlaD}).

A \emph {CW resolution} of a simplicial 
\Pa{\F} \ $\Ad$ \ is a CW complex \ $\Gd\in s\PAlg{\F}$, \ with CW basis \ 
$(\bar{G}_{n})_{n=0}^{\infty}$ \ such that each \ $\bar{G}_{n}$ \ is a free
\Pa{\F}, together with a weak equivalence \ $\phi:\Gd\to\Ad$.
\end{defn}

\begin{defn}\label{dcwr}\stepcounter{subsection}
In the situation of \S \ref{armc}, \ a simplicial object \ $\Rd\in s\C$ \ is 
called a \emph{CW resolution} of \ $\Xd\in s\C$ \ if \ $\Rd$ \ is a CW 
complex with each \ $\bar{R}_{n}$ \ in $\F$, up to homotopy (so in 
particular \ $\Rd$ \ is indeed cofibrant), equipped with a weak equivalence \ 
$f:\Rd\to\Xd$. 
\end{defn}

\begin{remark}\label{rcw}\stepcounter{subsection}
It is easy to see that one can inductively construct a CW resolution for 
every  simplicial \Pa{\F} \ $\Ad$, \ since in order for \ $\phi:\Gd\to\Ad$ \ 
to be a weak equivalence it is necessary and sufficient that \ $Z_{n}\phi$ \ 
take \ $Z_{n}\Gd$ \ onto a set of representatives of \ 
$\pi_{n}\Ad$ \ in \ $Z_{n}\Ad$, \ and \ 
the attaching map \ $\bar{d}^{n}_{0}$ \ map \ $\bar{G}_{n}$ \ onto 
a set of representatives for \ $\Ker(\pi_{n}\phi)$ \ in \ $Z_{n-1}\Gd$. \ 
Thus we can let \ $\bar{G}_{n}$ \ 
be the free \Pa{\F}\ (\S \ref{dfpa}) generated by union of the underlying sets 
of \ $Z_{n}\Ad$ \ and \ $\Ker(Z_{n-1}f)$, \ say. 
\end{remark}

The ``topological'' version of this requires a little more care. In particular,
\cite[Remark 3.16]{BlaCW} implies that not every free simplicial \Pa{\F} \ 
$\Ad$ \ is \emph{realizable} in the sense that there is a \ $\Rd\in s\C$ \ 
with \ $\pif\Rd\cong\Ad$. \ In order to see what can be said on this context,
assume given a fibrant and cofibrant simplicial object \ $\Pd$ \ with an 
augmentation \ $\varepsilon:P_{0}\to X$. \ For each \ $\alpha\in\hF$, \ 
consider the long exact sequence 
%
\setcounter{equation}{\value{thm}}\stepcounter{subsection}
\begin{equation}\label{ezero}
\dotsc \pi_{\alpha+1}C_{m}\Pd \xra{\ddi{m}} \pi_{\alpha+1}Z_{m-1}\Pd
\xra{\partial_{m-1}} \pia Z_{m}\Pd \xra{\jj{m}} \pia C_{m}\Pd 
\dotsc
\end{equation}
\setcounter{thm}{\value{equation}}
\noindent for the fibration \ $\dd{m}$, \ where \ $Z_{0}\Pd\DEF P_{0}$. \ 
By definition, \ $\Pd\to X$ \ is a resolution if and only if \ 
$\pi_{i}\pif\Pd=0$ \ for each \ $i\geq 0$, \  where the homotopy groups are
understod in the augmented sense \ -- \ that is, \ 
$\pi_{0}\pif\Pd\DEF \Ker(\ddi{0}:C_{0}\pif\Pd\to Z_{-1}\pif\Pd)
/\Image(\ddi{1}:C_{1}\pif\Pd\to Z_{0}\pif\Pd)$. \ 
The key technical fact we shall need in this context is contained in the following
%
%
\begin{lemma}\label{ltwo}\stepcounter{subsection}
An fibrant and cofibrant \ $\Pd\in s\C$ \ with an augmentation \ $\Pd\to X$ \ 
is a resolution of $X$ if and only if for each \ $m> 0$:

\begin{enumerate}
\renewcommand{\labelenumi}{(\alph{enumi})~}
\item There is a short exact sequence \ 
$0\to \Image(\partial_{m-1})\hra \pif Z_{m}\Pd 
\xra{\jj{m}} Z_{m}\pif\Pd \to 0$, \ and
\item $\partial_{m}\rest{\Image(\partial_{m-1})}$ \ 
is one-to-one, and surjects onto \ $\Image(\partial_{m})$, \ and \ 
$\Image\partial_{0}\cong\Omega\pif X$.
\end{enumerate}
\end{lemma}

Note that since \ $\partial_{m}$ \ shifts degrees by one, (a) and (b) together
imply that \ $\Image(\partial_{m})\cong \Omega^{m+1}\pif X$ \ for each $m$.

\begin{proof}
For any \ $\Pd$, \ the inclusion \ $j_{m}: Z_{m}\Pd\to C_{m}\Pd$ \ induces a map 
of \Pa{\F}s \ $\jj{m}:\pif Z_{m}\Pd\to \pif C_{m}\Pd\cong C_{m}\pif\Pd$ \ 
(see Lemma \ref{lone}), which factors through \ $Z_{m}\pif\Pd$. \ Denote 
the boundary map for the chain complex \ $C_{\ast}\pif\Pd$ \ (which computes \ 
$\pis\pif\Pd$) \ by \ $D_{m}\DEF\jj{m-1}\circ\ddi{m}$.

If \ $\Pd\to X$ \ is a resolution, we must have \ 
$\Image(\jj{m}\circ\ddi{m+1})=\Image(D_{m+1})=\Ker(D_{m})$ \ for each \ 
$m\geq 0$, \ so in particular \ $\jj{m}$ \ maps onto \ $Z_{m-1}\pif\Pd$. \  
Moreover, since \ $\pif C_{1}\Pd\to\pif P_{0}\epic\pif X\to 0$ \ is exact, \ 
$\Image\partial_{0}\cong\Omega\pif X$ \ and so if we assume by induction that
(b) holds for \ $m-1$, \ we see that \ $\Ker\jj{m}=\Image\partial_{m-1}$ \ is
isomorphic to \ $\Omega^{m}\pif X$, \ which proves (a). Moreover, if \ 
$0\neq\gamma\in\Ker\partial_{m}=\Image\ddi{m+1}$, \ and \ 
$\gamma\in\Image\partial_{m-1}=\Ker\jj{m}$, \ then we have \ 
$\beta\in\pif C_{m+1}\Pd$ \ with \ $\ddi{m+1}(\beta)=\gamma\neq 0$ \ but \ 
$D_{m+1}(\beta)=0$ \ -- \ contradicting (a) for \ $m+1$. \ Finally, if \ 
$\jj{m}(\gamma)\neq 0$, \ there is a \ $\beta\in\pif C_{m+1}\Pd$ \ with \ 
$D_{m}(\beta)=\jj{m}(\gamma)$, \ by the acyclicity of \ $\pif\Pd$, \ so \ 
$\gamma-\ddi{m+1}(\beta)\in\Ker\jj{m}=\Image\partial_{m-1}$, \ and \ 
$\partial_{m}(\gamma-\ddi{m+1}(\beta))=\partial_{m}(\gamma)$, \ which 
proves (b) for $m$. \ The identification of \ $\Image\partial_{0}$ \ is 
immediate from \eqref{ezero}.

Conversely, if (a) and (b) are satisfied for all $m$, for any element in \ 
$\zeta\in Z_{m}\pif\Pd$, \ we have \ $\zeta=\jj{m}(\gamma)$ \ for some \ 
$\gamma\in\pif Z_{m}\Pd$. \ Thus there is a \ $\theta\in\pif Z_{m-1}\Pd$ \ 
with \ $\partial_{m}(\partial_{m-1}(\theta))=\partial_{m}(\gamma)$, \ by (b), 
so \ $\gamma\cdot\partial_{m-1}(\theta)^{-1}$ \ is in \ 
$\Ker\partial_{m}=\Image\ddi{m+1}$; \ thus \ 
$\jj{m}(\gamma\cdot\partial_{m-1}(\theta)^{-1})=\zeta$ \ bounds, and \ 
$\pif\Pd$ \ is acyclic. 
\end{proof}

It should be pointed out that the fundamental short exact sequence \ 
%
\setcounter{equation}{\value{thm}}\stepcounter{subsection}
\begin{equation}\label{eeight}
0\to \Omega^{m}\pif X\cong \Image(\partial_{m-1})\hra \pif Z_{m}\Pd 
\xepic{\jj{m}} Z_{m}\pif\Pd \to 0
\end{equation}
\setcounter{thm}{\value{equation}}
\noindent for a resolution \ $\Pd$ \  is actually \emph{split}, as 
a sequence of graded groups, because \ 
$\jj{m}\rest{\Image\ddi{m+1}}=\jj{m}\rest{\Ker\partial_{m}}$ \ is one-to-one,
by (b), and surjects onto \ $Z_{m}\pif\Pd$ \ by the acyclicity. 
However, \ $\Image\ddi{m+1}=\Ker\partial_{m}$ \ need not be a sub-\Pa{\F}\ of \ 
$\pif Z_{m}\Pd$, \ since \ $\partial_{m}$ \ is not a morphism of \Pa{\F}s\vsm .

With the aid of Lemma \ref{ltwo} we can now show:
%
%
\begin{prop}\label{ptwo}\stepcounter{subsection}
Under the assumptions of \S \ref{armc} and \ref{agp}, any \ $X\in\C$ has a
CW resolution \ $\Rd\in s\C$.
\end{prop}

\begin{proof}
Let \ $\Qd\in s\C$ \ be the functorial resolution of \S \ref{sfres}; we may 
assume that the augmentation \ $\varepsilon^{Q}:Q_{0}\to X$ \ is a 
fibration.

We start off by choosing a set \ $\Ts^{0}\subseteq \pif Q_{0}$ \ of 
\Pa{\F}\ generators (\S \ref{dfpa}), such that if we let \ $R'_{0}\DEF
\coprod_{\alpha\in\hF}\coprod_{\beta\in T^{0}_{\alpha}}\ \Ma_{\beta}$, \ 
then \ $\varepsilon_{\#}^{Q}$ \ maps the free \Pa{\F} \ 
$\pif R'_{0}\subset \pif Q{_0}$ \ onto \ $\pif X$. \ We may assume \ 
$\Ts^{0}$ \ is \emph{minimal}, in the sense that no sub-graded set
generates a free \Pa{\F}\ surjecting onto \ $\pif X$ \ -- \ so that \ \
$\varepsilon_{\#}^{Q}(\beta)\neq 0$ \ for all \ $\beta\in \Ts^{0}$.

The inclusion \ $\phi:\pif R'_{0}\hra\pif Q{_0}$ \ defines a map \ 
$f'_{0}:R'_{0}\to Q_{0}$ \ with \ $(f'_{0})_{\#}=\phi$, \ and we let \ 
$\varepsilon^{R'}\DEF\varepsilon^{Q}\circ f'_{0}$; \ factoring \ 
$\varepsilon^{R'}$ \ by \ref{smc}(3) as \ 
$R'_{0}\xra{i}R_{0}\xra{\varepsilon^{R}}X$ \ and usng the LLP for $i$ and \ 
$\varepsilon^{Q}$ \ yields \ $f_{0}: R_{0}\to Q_{0}$ \ commuting with \ 
$\varepsilon$\vsm .
 
Now assume by induction that we have constructed a fibrant and 
cofibrant \ $\Rd$ \ through simplicial dimension \ $n-1\geq 0$, \ together 
with a map \ $\tr{n-1}f:\tr{n-1}\Rd\to\tr{n-1}\Qd$ \ which induces an 
embedding of \Pa{\F}s \ $(\tr{n-1}f)_{\#}$. \ 
We assume that \ $\Rd$ \ satisfies (a) and (b) of Lemma \ref{ltwo} for \ 
$0<m<n$ \ (and of course \ $\Qd$ \ satisfies them for all \ $m>0$). \ 
If we map the short exact sequence (a) for \ $\Rd$ \ to the
corresponding sequence for \ $\Qd$ \ by \ $f_{\ast}$, \ we see that \ 
$Z_{n-1}(f_{\#})=Z_{n-1}\phi:Z_{n-1}\pif\Rd\to Z_{n-1}\pif\Qd$ \ is 
one-to-one, so \ $(Z_{n-1}f)_{\#}:\pif Z_{n-1}\Rd\to\pif Z_{n-1}\Qd$ \ is, too.

Any non-zero element in \ $Z_{n-1}\pia\Rd$ \ is represented by \ 
$\gamma\in \pif Z_{n-1}\Rd$, \ by \eqref{eeight} for \ $R_{n-1}$. \ Let \ 
$g:\Ma\to Z_{n-1}\Qd$ \ represent \ $f_{\#}\gamma\in\pia Z_{n-1}\Qd$, \ 
with \ $\Ma_{(g)}$ \ the corresponding coproduct summand of \ 
$Q_{n}=LQ_{n-1}$ \ in \eqref{etwo}, \ with \ $i_{(g)}:\Ma_{(g)}\to Q_{n}$ \ 
the inclusion. \ Then \ $d_{i}\circ i_{(g)}=i_{(d_{i-1}g)}$ \ for \ 
$1\leq i\leq n$ \ (in the same notation) and \ $d_{0}\circ i_{(g)}=g$, \ 
by \S \ref{sfres}. Thus the \Pa{\F}\ generator \ 
$\lra{i_{(g)}}\in\pia Q_{n}$ \ is in \ $\C_{n}\pif\Qd$, \ and \ 
$\ddi{n}\lra{i_{(g)}}=f_{\#}\gamma$.  

Thus if we choose a set \ $\Ts^{n}$ \ of \Pa{\F} generators for \ 
$Z_{n-1}\pif \Rd$ \ and set
%
\setcounter{equation}{\value{thm}}\stepcounter{subsection}
\begin{equation}\label{eten}
\bar{R}_{n}\DEF\ \coprod_{\alpha\in\hF}\ \coprod_{\beta\in T^{n}_{\alpha}}\ 
\Ma_{(\beta)},
\end{equation}
\setcounter{thm}{\value{equation}}
\noindent we have maps \ $\bar{f}_{n}:\bar{R}_{n}\to C_{n}\Qd$ \ and \ 
$\bar{d}_{0}:\bar{R}_{n}\to Z_{n-1}\Rd$ \ such that \ 
$\jj{n-1}\circ\ddi{Q}\circ (\bar{f}_{n})_{\#}= 
\jj{n-1}\circ (Z_{n-1}f)_{\#}\circ(\bar{d}_{0})_{\#}$. \ Now \eqref{eeight}
implies that \ $\jj{n-1}$ \ is one-to-one on \ $\Image\bd$, \ so \ 
$\ddi{Q}\circ (\bar{f}_{n})_{\#}=(Z_{n-1}f)_{\#}\circ(\bar{d}_{0})_{\#}$. \ 
Because \ $\ddi{Q}$ \ is a fibration and \ $\pif\bar{R}_{n}$ \ is free, this
implies that one can choose \ $\bar{f}_{n}$ \ so that \ 
$\dd{Q}\circ \bar{f}_{n}= Z_{n-1}f\circ \bar{d}_{0}$. \ Since \ 
$L_{n}f:L_{n}\Rd\to L_{n}\Qd$ \ exists by the induction hypothesis, \ 
one can define \ $f_{n}:R_{n}\simeq L_{n}\Rd\amalg \bar{R}_{n}\to Q_{n}$ \ 
extending \ $\tr{n-1}f$ \ to \ 
$\tr{n}f:\tr{n}\Rd\to\tr{n}\Qd$, \ with \ 
$\delta^{R}_{n}:R_{n}\to M_{n}\Rd$ \ a fibration. \ Since \ $\pi_{i}\pif\Pd=0$ \ 
then holds for \ $i\leq n-1$, \ \eqref{eten} and \eqref{eeight} hold for \ $m=n$.
\end{proof}

\begin{remark}\label{rres}\stepcounter{subsection}
We have actually proved a little more: given any \emph{minimal} simplicial
CW resolution of \Pa{\F}'s \ $\Ad\to\pif X$ \ (\S \ref{dcw}) of a realizable
\Pa{\F}, \ one can find a CW resolution \ $\Rd\to X$ \ realizing it: 
that is, \ $\pif\Rd\cong\Ad$. \ (Minimality here is understood to mean that
we allow no unnecessary \Pa{\F}\ generators in each \ $\bar{A}_{n}$, \ beyond
those needed to map onto \ $Z_{n-1}\Ad$.)

By a more careful analysis, as in \cite[Thm.\ 3.19]{BlaCW}, one could in fact
show that \emph{any} CW resolution of \ $\pif X$ \ is realizable. However, 
this will follow from Corollary \ref{cone} below.
\end{remark}
%
%
\sect{Postnikov systems and the fundamental group action}
\label{cpf}

We now describe Postnikov systems for simplicial objects in the resolution model 
category, and the fundamental group action on them. 

\begin{defn}\label{dpt}\stepcounter{subsection}
If \ $\C$ \ is a category satisfying the assumptions of \S \ref{armc}, a 
\emph{Postnikov system} for an object \ $\Yd\in s\C$ \ is a sequence of objects \ 
$P_{n}\Xd\in s\C$, \  together with maps \ 
$\varphi^{n}:\Xd\to P_{n}\Xd$ \ and \ $p^{n}:P_{n+1}\Xd\to P_{n}\Xd$ \ (for \ 
$n\geq 0$), \ such that \ $\pin{k}{p^{n}}$ \ and \ $\pin{k}{\varphi^{n}}$ \ are 
isomorphisms for all \ $k\leq n$, \ and \ $\pin{k}{P_{n}\Xd}=0$ \ for \ 
$k\geq n+1$
\end{defn}

\begin{remark}\label{rpt}\stepcounter{subsection}
In general, such Postnikov towers may be constructed for fibrant \ $\Xd$ \ using 
a variant of the standard construction for simplicial sets 
(cf.\ \cite[\S 8]{MayS}) due to Dwyer and Kan in \cite[\S 1.2]{DKaO}, \ and 
for arbitrary \ $\Xd$ \ by using a fibrant approximation. 

Note that if \ $\Qd\in s\C$ \ is a resolution of some \ $X\in\C$ \ \
(see \S \ref{sres}), \ then by \eqref{eseven} \ 
$\pin{i}{P_{n}\Qd}\cong \Omega^{i}\pif X$ \ for \ $n\geq i\geq 0$, \ and \ 
$\pin{i}{P_{n}\Qd}=0$ \ for \ $i>n$; \ so \eqref{esix} implies that \ 
%
\setcounter{equation}{\value{thm}}\stepcounter{subsection}
\begin{equation}\label{eeleven}
\pi_{i}\pif P_{n}\Qd\cong\begin{cases}
\pif X&\text{for \ }i=0\\
\Omega^{n+1}\pif X&\text{for \ }i=n+2,\\
0 & \text{otherwise.}
\end{cases}
\end{equation}
\setcounter{thm}{\value{equation}}
\end{remark}

\subsection{Postnikov towers for resolutions}
\label{sptr}\stepcounter{thm}

\noindent It is actually easier to construct a cofibrant version of the 
Postnikov tower for a resolution than it is to construct the resolution itself:
Given a CW resolution \ $\Qd$ \ of an object \ $X\in\C$, \ (constructed as
in Proposition \ref{ptwo}), with CW basis \ $(\bar{Q}_{k})_{k=0}^{\infty}$, \ 
we construct a CW cofibrant approximation \ $\Yd\to \Qd\q{n}$ \ as follows.

Let \ $\Js\DEF\pif X$, \ and choose some \ $G\in ho\F$ \ (i.e., \ 
$G\simeq \coprod_{\alpha\in\hF}\ \coprod_{T_{\alpha}}\ \Ma$) \ having a 
surjection of \Pa{\F}s \ $\phi:\pif G\epic\Omega^{n+1}\Js$. \ Set \ 
$\bar{Y}_{n+2}\DEF \bar{Q}_{n+2}\amalg G$, \ with \ 
$(d_{0}\rest{G})_{\#}=\phi$, \ mapping onto \ 
$\Omega^{n+1}\Js\cong\Image(\partial_{n})\hra\pif Z_{n+1}\Qd=\pif Z_{n+1}\Yd$ \
(see \eqref{eeight}). \ This defines \ 
$Y'_{n+2}\DEF \bar{Y}_{n+1}\amalg L_{n+2}\Yd\xra{\delta_{n+2}} M_{n+2}\Yd$, \ 
which we then change into a fibration. Since \ 
$\ddi{n+2}:\pif C_{n+2}\Yd\epic \pif Z_{n+1}\Yd$ \ is surjective, we may
assume by induction on \ $k\geq n+2$ \ that
%
\setcounter{equation}{\value{thm}}\stepcounter{subsection}
\begin{equation}\label{etwelve}
(j_{k})_{\ast}:\pif Z_{k}\Yd\xra{\cong}Z_{k}\pif\Yd\hsp \text{and}\hsp
\partial_{k-1}=0,
\end{equation}
\setcounter{thm}{\value{equation}}

\noindent and thus we may choose \ 
$\bar{Y}_{k+1}\in ho\F$ \ with \ 
$\bar{d}_{0}:\pif\bar{Y}_{k+1}\epic\pif Z_{k}\pif\Yd$, \ and see that
\eqref{etwelve} holds for \ $k+1$ \ by \eqref{eeight}.

Note that \ $\Yd\simeq \Qd\q{n}$ \ is constructed by ``attaching cells'' to \ 
$\Qd$, \ as in the traditional method for ``killing homotopy groups'' 
(cf.\ \cite[\S 17]{GrayH}), so we have a natural \emph{embedding} \ 
$\rho:\Qd\hra \Yd$, \ rather than a fibration.
In fact, it is helpful to think of \ $P_{n}\Xd$ \ as a homotopy-invariant version
of the \ $(n+1)$-skeleton of \ $\Xd$: \ starting with \ $\tr{n+1}\Xd$, \ one 
completes it to a full simplicial object by a functorial construction which 
(unlike the skeleton) depends only on the homotopy type of \ $\Xd$.

\subsection{$\Pi$-algebras and the fundamental group}
\label{spf}\stepcounter{thm}

Under our assumptions, the category \ $\C=\PAlg{\F}$ \ is a \emph{CUGA}, or 
category of universal graded algebras (see \cite[\S 2.1]{BStG} and 
\cite[V,\S 6]{MacC}), so that \ $s\C$, \ the category of simplicial
\Pa{\F}s, has a model category structure defined by Quillen 
(see \cite[II, \S 4]{QuH}). \ Equivalently, one could take the resolution 
model category on \ $s\C$, \ starting with the trivial model category 
structure on \ $\PAlg{\F}$, \ and letting \ $\F_{\PAlg{\F}}$ \ be the 
subcategory of all free \Pa{\F}s \ -- \ as in \S \ref{ermc}. \ One thus has a 
concept of ``spheres'' in \ $s\PAlg{\F}$ \ -- \ namely, \ 
$\pif\Sigma^{n}\Ma$, \ for \ $\alpha\in\hF$ \ (cf.\ \S \ref{dsfs}) \ -- \ 
and \ $(\pi_{n}\Ad)_{\alpha}\cong [\Sigma^{n}\Ma, \Ad]_{s\PAlg{\F}}$ \ for 
any \ simplicial \Pa{\F}\ \ $\Ad$. \ Thus if we take homotopy classes of maps 
between (coproducts of) these spheres as the primary homotopy operations 
(see \cite[XI, \S 1]{GWhE}), we can endow the homotopy groups \ 
$\pis\Ad=(\pi_{i}\Ad)_{i=0}^{\infty}$ \ of \ $\Ad$ \ with an additional 
structure: that of a \ ($\PAlg{\F}$)-\Pa{}, \ in the (somewhat unfortunate, 
in this case) terminology of \cite[\S 3.2]{BStG}. By definition, this structure 
is a homotopy invariant of \ $\Ad$. \ 

In our situation, however, because we are dealing with Postnikov sections, by
\eqref{eeleven} we only need the very simplest part of that structure \ -- \ 
namely, the action of the fundamental group \ $\pi_{0}\Ad$ \ on each of the 
higher homotopy groups \ $\pi_{n}\Ad$. \ 

Observe that because $\C$ has an underlying group structure, by assumption 
\ref{agp}, the indexing of the homotopy groups of an object in \ $s\C$ \ should 
be shifted by one compared with the usual indexing in \ $\Ta$, \ so that \ 
$\pi_{0}\Ad$ \ is indeed the fundamental group, and in fact the action we refer
to is a straightforward generalization of the usual action of the fundamental 
group of a simplicial group (or topological space) on the higher homotopy groups.

\subsection{$\Js$-modules and \ $\Js$-algebras}
\label{sma}\stepcounter{thm}

We shall be interested in an algebraic description of this action: that is, we 
would like a category of universal algebras which model this action, in the same 
sense that \Pa{}s model the action of all the primary homotopy operations on the 
homotopy groups of a space. Just as in the case of ordinary \Pa{}s, 
the action in question is determined by the homotopy classes of maps of 
simplicial \Pa{\F}s.

Thus we are led to consider two distinct ``varieties of algebras'',
in the terminology of \cite[V, \S 6]{MacC}): one modeled on the homotopy classes
of maps, and one on the actual maps.

\begin{defn}\label{dmod}\stepcounter{subsection}
Given a \Pa{\F}\ \ $\Js$, \ let \ $\RM{\Js}$ \ denote the category of universal
algebras whose operations are in one-to-one correspondence with homotopy
classes of maps \ 
$\pif\Sigma^{n}\Ma\to\pif(\Sigma^{n}\M{\alpha'}\amalg\Sigma^{0}\M{\alpha''})$, \ 
and whose universal relations correspond to the relations holding among these 
homotopy class in \ $ho(s\C)$. \ These model \ $\pi_{n}\Ad$, \ with the action
of \ $\pi_{0}\Ad$, \ for \ $\Ad\in s\PAlg{\F}$.

An object \ $\Ks\in\RM{\Js}$ \ is itself a \Pa{\F}, \ equipped with an action 
of an operation \ $\lambda:J_{\alpha''}\times K_{\alpha'}\to K_{\alpha}$ \ 
for each \ $\lambda\in
[\pif\Sigma^{n}\Ma,\pif(\Sigma^{n}\M{\alpha'}\amalg\Sigma^{0}\M{\alpha''})]$. \ 
Such a \ $\Ks$  will be called a \ $\Js$-\emph{module}, even though in
general the category of such objects, which we shall denote by \ $\RM{\Js}$, \ 
need not be abelian (and it could depend on $n$). \ However, in the cases that
interest us, \ $\RM{\Js}$ \ will be abelian, and will not depend on \ $n>0$.
\end{defn}

\begin{defn}\label{dalg}\stepcounter{subsection}
Given a \Pa{\F}\ \ $\Js$, \ let \ $\Alg{\Js}$ \ denote the category of universal
algebras whose operations are in one-to-one correspondence with 
\emph{actual maps} \ 
$\pif\Sigma^{n}\Ma\to\pif(\Sigma^{n}\M{\alpha'}\amalg\Sigma^{0}\M{\alpha''})$ \ 
as above, and whose universal relations correspond to the relations holding 
among these maps in \ $s\C$. \ The objects in \ $\Alg{\Js}$, \ which are again
\Pa{\F}s with additional structure, will be called \ $\Js$-\emph{algebras}.
\end{defn}

The category \ $\Alg{\Js}$ \ is generally very complicated; it is not abelian, 
and we cannot expect to know much about it, even for \ $\C=\G$, \ say. In 
particular, one may well have a different category for each \ $n>0$ \ 
(although we surpress the dependence on $n$ to avoid excessive notation). 
Note, however, that maps \ $\ell:\pif\Sigma^{n}\Ma\to\Ad$, \ for any simplicial 
\Pa{\F}\ \ $\Ad$, \ correspond to elements in \ $Z_{n}\Ad$, \ so that the \ 
$A_{0}$-algebra structure on \ $A_{n}$ \ restricts to an action of
of \ $Z_{0}\Ad=A_{0}$ \ on \ $Z_{n}\Ad$. 

\begin{remark}\label{rloop}\stepcounter{subsection}
Let \ $\Qd$ \ be a resolution (in \ $s\C$) \ of some object \ $X\in\C$, \ 
with \ $\Js\DEF\pif X$, \ and \ 
$\Yd\simeq P_{n}\Qd$ \ its $n$-th Postnikov approximation. Then we have an
action of \ $\pi_{0}\pif\Yd\cong\Js$ \ on \ 
$\pi_{n+2}\pif\Yd\cong\Omega^{n+1}\Js$ \ which is a homotopy invariant of \ 
$\Yd$, \ and thus in turn of \ $\Qd$, \ so of $X$. \ It is not clear 
on the face of it whether the \ $\Js$-module \ $\Omega^{n}\Js$ \ depends only 
on \ $\Js$, \ though we shall see (in \S \ref{sco} below) this holds for \ 
$n=1$, \ and hope to show in \cite{BGoeC} that in fact this holds for all $n$. \ 
In any case it is describable purely in terms of the primary \Pa{\F}-structure 
of \ $\Js$. 

In general, for any simplicial object \ $\Xd\in s\C$, \ there is an action of \ 
$\pin{0}{\Xd}\cong\pi_{0}\pif\Xd$ \ on the higher \Pa{\F}s \ $\pin{n}{\Xd}$, \ 
defined similarly via homotopy classes of maps \ 
$[\SSp{n},\SSp{0}\amalg \SSp{n}]_{s\C}$ \ (see \S \ref{dsfs}); but there is no 
reason why this should define the same category of ``$\pin{0}{\Xd}$-modules'' 
as that defined above. Thus we do not know \eqref{esix}to be a long exact 
sequence of $\pin{0}{\Xd}$-modules. However, in our case, when \ $\Xd=\Qd$ \ is
a resolution, the isomorphism of (abelian) \Pa{\F}s \ 
$\pi_{n+2}\pif\Yd\cong\Omega^{n+1}\Js$ \ is \emph{defined} inductively by means 
of the connecting homomorphism of \eqref{esix}, and this yields the \ $\Js$-module
structure on \ $\Omega^{n}\Js$.
\end{remark}

\begin{assume}\label{aab}\stepcounter{subsection}
Under mild assumptions on the category $\C$ one may show that for any \ 
$\Ad\in s\PAlg{\F}$ \ and \ $n\geq 1$, \ the \Pa{\F}\ \ $\pi_{n}\Ad$ \ is 
abelian (see \cite[Lemma 5.2.1]{BStG}). 

However, we shall need to assume more than this: namely, that \ $\RM{\Js}$ \ 
as defined above is in fact an abelian category. We also assume that when \ 
$\Ad$ \ is a simplicial \Pa{\F}, \ the action of \ $\pi_{0}\Ad$ \ on each \ 
$\pi_{n}\Ad$ \ is induced by an action of \ $A_{0}$ \ on \ $A_{n}$, \ and if \ 
$\Ad=\pif\Qd$, \ then this in turn is induced by an action of \ $Q_{0}$ \ on \ 
$Q_{n}$. \ Moreover, \ $Z_{n}\Ad$ \ and \ $C_{n}\Ad$ \ are sub-$A_{0}$-algebras 
of \ $A_{n}$, \ and \ $\bd$ \ is a homomorphism of \ $A_{0}$-algebras.
\end{assume}
%
%
\begin{prop}\label{pthree}\stepcounter{subsection}
These assumptions are satisfied for the categories listed in \S \ref{ragp}.
\end{prop}

\begin{proof}
As we shall see, all the categories in question are essentially special cases 
of the first:

\begin{enumerate}
\renewcommand{\labelenumi}{(\Roman{enumi})~}
\item When \ $\C=\G$, \ the fundamental group action has an explicit 
description as follows:

We define the \emph{generalised Samelson product} of two elements \ 
$x\in X_{p,k}$ \ $y\in X_{q,\ell}$ \ (where, as in \S \ref{dso}, $p$ is the 
``external'' dimension, $k$ the ``internal'' dimension in a a bisimplicial 
group \ $\Xdd\in s\G$) \ to be the element \ $\llra{x,y}\in X_{p+q,k+\ell}$ \ 
%
\setcounter{equation}{\value{thm}}\stepcounter{subsection}
\begin{equation}\label{ethirteen}
\llra{x,y}\DEF
\prod_{(\sigma,\rho)\in S_{p,q}}\left(\prod_{(\varphi,\psi)\in S_{k,\ell}} 
(s^{ext}_{\rho_{q}}\dotsc s^{ext}_{\rho_{1}}
s^{int}_{\psi_{\ell}}\dotsc s^{int}_{\psi_{1}} x, 
s^{ext}_{\sigma_{q}}\dotsc s^{ext}_{\sigma_{1}}
s^{int}_{\varphi_{\ell}}\dotsc s^{int}_{\varphi_{1}} y
)^{\varepsilon(\varphi)}\right)^{\varepsilon(\sigma)}.
\end{equation}
\setcounter{thm}{\value{equation}}

Here \ $S_{p,q}$ \ is the set of all \ $(p,q)$-shuffles  \ -- \ 
that is, partitions of \ $\{0,1,\dotsc,p+q-1\}$ \ into disjoint sets \
$\sigma_{1}<\sigma_{2}<\dotsb <\sigma_{p}$, \ 
$\rho_{1}<\rho_{2}<\dotsb <\rho_{q}$ \ -- \ and \ $\varepsilon(\sigma)$ \ is 
the sign of the permutation corresponding to \ $(\sigma,\rho)$ \ 
(see \cite[VIII, \S 8]{MacH}); \ $S_{p,q}$ \ is ordered by the reverse 
lexicographical ordering in $\sigma$. \ $(a,b)$ \ denotes the commutator \ 
$a\cdot b\cdot a^{-1}\cdot b^{-1}$ \ (where \ $\cdot$ \ is the group 
operation).
When \ $p=q=0$, \ $\llra{x,y}$ \ is just the usual Samelson product \ 
$\lra{x,y}$ \ in \ $X_{0,\bullet}\in\G$ \ (cf.\ \cite[\S 11.11]{CurtS}).

We are mainly interested here in the case \ 
$p=0$, \ so \ $\llra{x,y}\DEF \lra{\hat{x},y}$ \ for \ 
$\hat{x}\DEF s_{q-1}\dotsb s_{0}x\in X_{q,k}$. \ It is sometimes
convenient to think of this as an ``action'' of $x$ on $y$, setting \ 
$t_{x}(y)\DEF \llra{x,y}\cdot y$ \ (cf.\ \cite[X, (7.4)]{GWhE}). 

The simplicial identities imply that if \ 
$d^{int}_{i}x= d^{int}_{i}y=0$ \ for all $i$, \ the same holds for \ 
$\lra{x,y}$, \ and if \ $x=d^{int}_{0}z$ \ for some \ 
$z\in C_{k+1}^{int}X_{p,\bullet}$, \ then \ $\lra{x,y}=d_{0}^{int}\lra{z,y}$, \ 
so that \ $\llra{\ ,\ }$ \ induces a well-defined operation \ 
$\llra{\ ,\ }:\pi^{int}_{k}X_{p,\bullet}\times\pi^{int}_{\ell}X_{q,\bullet}\to
\pi^{int}_{k+\ell}X_{p+q,\bullet}$, \ which is defined for any simplicial \Pa{}\ 
$\Ads$, \ with \ $\alpha\in A_{p\ast}$ \ and \ $\beta\in A_{q\ast}$, \ by:
%
\setcounter{equation}{\value{thm}}\stepcounter{subsection}
\begin{equation}\label{efourteen}
\llra{\alpha,\beta}\DEF\prod_{(\sigma,\rho)\in S_{p,q}}\ 
\lra{s_{\rho_{q}}\dotsc s_{\rho_{1}}\alpha,
s_{\sigma_{q}}\dotsc s_{\sigma_{1}}\beta}^{\varepsilon(\sigma)}\in A_{p+q\ast}.
\end{equation}
\setcounter{thm}{\value{equation}}

Again when \ $p=0$ \ 
we write \ $\tau_{\alpha}(\beta)\DEF \llra{\alpha,\beta}\cdot \beta$, \ so that \ 
$\tau_{\alpha}:A_{q\ast}\to A_{q\ast}$ \ is a group homomorphism in each 
degree (if \ $\alpha\in Z_{p}\Ads$, \ $\beta\in Z_{q}\Ads$, \ 
then \ $\llra{\alpha,\beta}\in Z_{p+q}\Ads$).

Now let \ $\Xdd\DEF \Sigma^{0}\bS{k}\amalg \Sigma^{n}\bS{\ell}$, \ 
(where \ $\bS{k}$ \ is the $k$-sphere for $\G$ \ -- \ \S \ref{ermc}) and let \ 
$\iota_{0,k}$ \ and $\iota_{n,\ell}$ \ be \Pa{}-generators for \ 
$\pi_{k}X_{0,\bullet}$ \ and \ 
$\pi_{\ell}\bS{\ell}\subseteq\pis X_{n,\bullet}$, \ respectively, so \ 
$\pis X_{n,\bullet}$ \ is generated by \ 
$\{\hat{\iota}_{0,k},\iota_{n,\ell}\}$. \ 
Since \ $d_{j}\iota_{n,\ell}=0$ \ ($0\leq j\leq n$), \ we have a short exact 
sequence of \Pa{}s
%
\setcounter{equation}{\value{thm}}\stepcounter{subsection}
\begin{equation}\label{eeighteen}
0\to Z_{n}\pis\Xdd\to\pis(\bS{k}\amalg\bS{\ell})\to \pis\bS{k}\to 0.
\end{equation}
\setcounter{thm}{\value{equation}}

When \ $k,\ell>0$, \ by \cite[Theorem A]{HilH} any element \ 
$x\in\pis X_{n,\bullet}\cong\pis(\bS{k}\amalg\bS{\ell})$ \ 
can be written as a sum of elements of the form \ 
$\zeta^{\#}\omega(\hat{\iota}_{0,k},\iota_{n,\ell})$ \ (where \ 
$\omega(x,y)=\lra{\dotsc\lra{x,y},\dotsc}$ \ is some iterated Samelson 
product), so $x$ can be obtained by means of the ``internal'' 
$\PAlg{}$\ operations from expressions of the form \ 
$\tau_{\alpha}(\iota_{n,\ell})$ \ (for \ $\alpha\in\pis X_{0,\bullet}$).

By passing to universal covers we have a similar description when \ $\ell>k=0$, \ 
since then any \ $x\in\pi_{j}X_{n,\bullet}$ \ ($j\geq 1$) \ can be written as a 
sum of elements of the form \ 
$\zeta^{\#}\omega(
\tau_{\alpha_{1}}(\iota_{n,\ell}),\dotsc,\tau_{\alpha_{r}}(\iota_{n,\ell})$ \ 
(for \ $\alpha_{i}\in\pis X_{0,\bullet}$), \ and any other \ 
$\alpha\in\pis X_{0,\bullet}$ \ acts on this by permuting the generators \ 
$\tau_{\alpha_{i}}(\iota_{n,\ell}$, \ so again \ $tau_{\alpha}(-)$ \ is a group
homomorphism. When \ $k=\ell=0$, \ we are reduced to the case \ $\C=\Gp$ \ 
(see (II) below).

When \ $k>0$ \ and \ $\ell=0$, \ let us write \ 
$\varphi_{\alpha}(\beta)\DEF\llra{\beta,\alpha}$ \ 
for \ $\alpha\in\pis X_{0,\bullet}$ \ and  \ $\beta\in\pi_{0}X_{n,\bullet}$, \ 
so that we are thinking of the usual (internal) action 
of the fundamental group \ $\pi_{0}X_{n,\bullet}$ \ as a function of \ $\beta$. \ 
This is \emph{not} a homomorphism, since we have \ 
$\varphi_{\alpha}(\beta\cdot\gamma)=
\varphi_{\alpha}(\beta)+\varphi_{\alpha}(\gamma)+
\llra{\beta,\llra{\gamma,\hat{\alpha}}}$ \ by \cite[III, (1.7) \& X, (7.4)]{GWhE}.

But \ $\llra{\alpha,\beta}$ \ is a cycle (i.e., in \ $Z_{n}\pis\Xdd$), \ by 
\eqref{eeighteen}, so \ $\llra{\llra{\alpha,\beta},\gamma}\sim 0$ \ in \ 
$\pi_{n}\pis\Xdd$ \ for any \ $\gamma\in\pi_{0}X_{n,\bullet}$ \ by 
\cite[5.2.1]{BStG}, \ which means that \ $\varphi_{\alpha}$ \ induces
a homomorphism on \ $\pi_{n}\pis\Xdd$.

In summary, an $\Js$-\emph{algebra} (\S \ref{dalg}), for any \ $\Js\in\PAlg{}$, \ 
is just a \Pa{}\ \ $\Ks$ \ together with an action of each \ $\alpha\in\Js$, \ 
which may be expressed in terms of the (degree-shifting) homomorphisms \ 
$\tau_{\alpha}$, \ or the functions \ $\varphi_{\alpha}$, \ respectively, 
satisfying whatever relations hold among these (and the internal \Pa{}\ 
operations) in \ $\pis X_{n,\bullet}$.

A \ $\Js$-\emph{module}, on the other hand, is an \emph{abelian} \Pa{}\ \ 
$\Ks$, \ together with homomorphisms \ $\tau_{\alpha}:\Ks\to\Ks$ \ or \ 
$\varphi_{\alpha}:\Ks\to\Ks$ \ for each \ $\alpha\in\Js$, \ satisfying the 
identities occuring in \ $\pi_{n}\pis\Xdd$. \ 

These identities could be described more or less explicitly in the category \ 
$\PAlg{}$, \ in terms of suitable Hopf invariants (cf.\ \cite[II, \S 3]{BauC}). 
Compare \cite[\S 3]{BauCC}). 
\item When \ $\C=\Gp$, \ $s\C$ \ models the homotopy theory of (connected) 
topological spaces, and \ $\RM{\Js}$, \ defined (as noted above) through the 
usual action of the fundamental group, is equivalent to the category of (left) 
modules over the group ring \ $\Z[\pi_{0}\Ad]$ \ (for \ $\Ad\in s\C\approx \G$).
\item When \ $\C=\Lie$, \ the situation is similar to \ $\C=\Gp$, \ with
Samelson products replaced by Lie brackets.
\item When \ $\C=d\LL\approx s\Lie$, \ one has a generalized Lie bracket defined 
for bisimplicial Lie algebras as in \eqref{ethirteen}, with commutators 
replaced by Lie brackets (see \cite[\S 2.6]{BlaHR}).
\item When \ $\C=\RM{R}$, \ $s\C$ \ is equivalent to the category of chain 
complexes over $R$, so there is no action of \ $\pi_{0}=H_{0}$ \ on the higher 
groups.
\end{enumerate}
\end{proof}

\begin{remark}\label{ruga}\stepcounter{subsection}
It is possible to write down general conditions on category of universal 
algebras (or CUGA) $\C$, defined in terms of operations and relations, which 
suffice to ensure that assumptions \ref{aab} hold: all one really needs is a
suitable Hilton-Milnor theorem in \ $s\C$ \ (see, e.g., \cite{GoeH}). However,
it seems simpler to state the conditions needed as above, and verify them 
directly in any particular case of interest.
\end{remark}
%
%
\sect{Cohomology of \Pa{\F}s}
\label{ccpa}

In this section we complete the description of the algebraic invariants 
used to distinguish homotopy types. To do so, we recall Quillen's definition 
of cohomology in a model category, in the context of \ $\PAlg{\F}$:

\begin{defn}\label{dhom}\stepcounter{subsection}
Let $\C$ be a model category with an \emph{abelianization} functor \ 
$\Ab:\C\to \C_{ab}$, \ where \ $\C_{ab}$ \ denotes of course the full category 
of abelian objects in $\C$; \ we shall usually write \ $X_{ab}$ \ for \ 
$\Ab(X)$ \ (see \S \ref{dapa}). \ In \cite[II, \S 5]{QuH} 
(or \cite[\S 2]{QuC}), Quillen defines the \emph{homology} of an object \ 
$X\in\C$ \ to be the total left derived functor \ $\mathbf{L}\Ab$ \ of \ 
$\Ab$, \ applied to $X$ \ (cf.\ \cite[I, \S 4]{QuH}). Likewise,
given an object \ $M\in\C_{ab}/X$, \ the \emph{cohomology of $X$ with 
coefficients in $M$} is \ 
$\mathbf{R}\Hom_{\C_{ab}/X}(X,M)\DEF\Hom_{\C_{ab}/X}(\mathbf{L}\Ab X,M)$. \ 
\end{defn}

\subsection{Quillen cohomology of \Pa{\F}s}
\label{sqcp}\stepcounter{thm}

When \ $\Js\in\C=\PAlg{\F}$, \ we have the model category structure defined in 
\S \ref{spf} above,  so we can choose a resolution \ $\Ad\to\Js$ \ 
in \ $s\PAlg{\F}$ \ as in \S \ref{sres}, and define the $i$-th homology 
group of \ $\Js$ \ to be the $i$-th homotopy group \ $\pi_{i}(\Ab\Ad)$ \ of the
$\hF$-graded simplicial abelian group \ $(\Ad)_{ab}$ \ -- \ i.e., of the
associated chain complex (cf.\ \cite[\S 1]{DoldH}).
One must verify, of course, that this definition is independent of the
choice of the resolution \ $\Ad\to \Js$.

Similarly, if \ $\Ks$ \ is an abelian \ $\Js$-algebra, then the $i$-th cohomology 
group of \ $\Js$ \ with coefficients in \ $\Ks$, \ written \ $H^{i}(\Js;\Ks)$, \ 
is that of the cochain complex corresponding to the cosimplicial 
$\hF$-graded abelian group \ $\Hom_{\Alg{\Js}}(\Ad,\Ks)$.

\begin{remark}\label{rrpa}\stepcounter{subsection}
Here \ $\Hom_{\Alg{\Js}}(A,B)$ \ is the group of \Pa{\F} homomorphisms which
respect the \ $\Js$-action; because we are mapping into an abelian object \ 
$\Ks$, \ $\Hom_{\Alg{\Js}}(\Ad,\Ks)\cong\Hom_{\Alg{\Js}}((\Ad)'_{ab},\Ks)$ \ 
(where \ $A'_{ab}$ \ denotes the abelianization of \ $A\in\Alg{\Js}$ \ as an 
$\Js$-algebra).

However, in the simplicial abelian \ $\Js$-algebra \ $(\Ad)'_{ab}$ \ we have a 
direct product decomposition \ 
$(A_{k})'_{ab}=(\hat{A}'_{k})_{ab}\oplus (L_{k}\Ad)'_{ab}$ \ for \ $k\geq 0$, \ 
where \ $(\hat{A}_{k})'_{ab}\DEF C_{k}(\Ad)'_{ab}$ \ is the the sub-abelian 
$\Js$-algebra of \ $(A_{k})'_{ab}$ \ generated by \ $(\bar{A}_{k})'_{ab}$ \ 
(cf.\ \cite[Cor.\ 22.2]{MayS}) \ -- \ and in fact \ 
$(\hat{d}_{0})'_{ab}:(\hat{A}_{n})'_{ab}\to(A_{n-1})'_{ab}$ \ factors through a
map \ $\hat{\partial}_{n}:(\hat{A}'_{n})_{ab}\to(\hat{A}_{n-1})'_{ab}$ \ 
(see \cite[p.\ 95(i)]{MayS}).

Thus the $n$-cochains split as:
$$
\Hom_{\Alg{\Js}}((A_{n})'_{ab},\Ks)\cong 
\Hom_{\Alg{\Js}}((\hat{A}_{n})'_{ab},\Ks)\oplus 
\Hom_{\Alg{\Js}}(L_{n}(\Ad)'_{ab},\Ks),
$$
\noindent so by \cite[X, \S 7.1]{BKaH} any cocycle representing a cohomology 
class in \ $H^{n}(\Js;\Ks)$ \ may be represented uniquely either by a map of 
abelian \ $A_{0}$-algebras \ $\hat{f}:(\hat{A}_{n})'_{ab}\to \Ks$, \ or by a 
map of \ $A_{0}$-algebras \ $f:A_{n}\to\Ks$.

Since \ $C_{n}\Ad$ \ contains the sub-$A_{0}$-algebra of \ $A_{n}$ \ generated 
by \ $\bar{A}_{n}$ \ (by assumption \ref{aab}), $f$ determines its restriction \ 
$f\rest{C_{n}\Ad}:C_{n}\Ad\to \Ks$, \ which determines \ $\hat{f}$, \ which
determines $f$ in turn. We have thus shown that \ $H^{\ast}(\Js;\Ks)$ \ may be
calculated as the cohomology of the (abelian) cochain complex \ 
$\Hom_{\Alg{A_{0}}}(\Cs\Ad,\Ks)$ \ (even though \ $\Cs\Ad$ \ is not in general
a homotopy invariant of \ $\Ad$, \ in non-abelian categories).
\end{remark}

\subsection{obstructions to existence of resolutions}
\label{soer}\stepcounter{thm}

Given an object \ $X\in\C$, \ and a (suitable) simplicial resolution \ 
$\Ad\to\Js$ \ of the \Pa{\F}\ $\Js\DEF\pif X$, \ we have seen in Section \
\ref{cmc} that one can construct a resolution \ $\Qd$ \ of $X$ (in the 
resolution model category \ $s\C$) \ realizing \ $\Ad$, \ in the sense that \ 
$\pif\Qd\cong\Ad$. \ It is thus natural to ask whether any simplicial 
\Pa{\F} \ -- \ or at least, any resolution \ $\Ad$ \ of an abstract \Pa{\F}\ \ 
$\Js$ \ -- \ is realizable in \ $s\C$. 

One approach to this question in the topological setting (i.e., for \ 
$\C=\G$), \ in terms of higher homotopy operations, was given in \cite{BlaHH}.
However, a glance at the proof of Proposition \ref{ptwo} shows that one can
instead consider obstructions to extending \ $\tr{n}\Qd$ \ to the next simplicial 
dimension. For a homotopy-invariant description, we 
state this in terms of successive Postnikov approximations to \ $\Qd$, \ since 
it is clear that, once we have constructed \ $\tr{n}\Qd$, \ it is always 
possible to obtain \ $\Yd\simeq \Qd\q{n-1}$ \ from it by successive choices of 
free objects \ $\bar{Y}_{k+1}\in ho\F$ \ ($k=n,\dotsc$) \ mapping to \ 
$Z_{k}\Yd$ \ by a \Pa{\F}\ surjection.

\subsection{constructing the obstruction}
\label{sco}\stepcounter{thm}

Assume given a CW resolution \ $\Ad\in s\PAlg{\F}$ \ of \ $\Js$, \ 
with CW basis \ $(\bar{A}_{n})_{n=0}^{\infty}$, \ and choose corresponding
free objects \ $\bar{Q}_{n}\in\F\subset\C$ \ with \ 
$\pif\bar{Q}_{n}\cong\bar{A}_{n}$. \ We begin the induction with \ 
$\tr{1}\Qd$, \ and thus \ $\Qd\q{0}$, \ constructed as in the proof of 
Proposition \ref{ptwo}. Note that to obtain \ $\tr{1}\Qd$ \ we do not in fact
need to know \ $X\in\C$ \ with \ $\pif X\cong\Js$ \ -- \ or even to know that 
such an object exists! This implies that the \ $\Js$-module structure on \ 
$\Omega\Js$ \ is uniquely determined.

In the inductive stage we assume given \ $\tr{n}\Qd$ \ (equivalently: \
$\Qd\q{n-1}$), \ satisfying \ref{ltwo}(a) and (b) for \ 
$0<m\leq n$. \ Our strategy is to try to attach \ $(n+1)$-dimensional ``cells'' 
to \ $\tr{n}\Qd$ \ in such a way as to guarantee acyclicity of the resulting \ 
$\tr{n+1}\Qd$ \ in one more simplicial dimension \ -- \ using Lemma \ref{ltwo}
above. The key to the construction of \ $\tr{n+1}\Qd$ \ from \ $\tr{n}\Qd$ \ 
thus lies in the extension of \ $A_{0}$-algebras \ \eqref{eeight} (for \ $\Qd$, \ 
rather than \ $\Pd$), \ in which the two ends are given to us. Observe that this
extension determines the \ $A_{0}$-algebra structure on \ $\Omega^{n}\Js$, \ if 
more than one is possible.

We want this extension to be ``trivial'' (that is, split as a semi-direct 
product of \ $A_{0}$-algebras), in order to be able to lift the given map of \ 
$A_{0}$-algebras \ $\bar{d}_{0}^{A}:\bar{A}_{n+1}\to Z_{n}\Ad$ \ to a map \ 
$\bar{d}_{0}^{Q}:\bar{Q}_{n+1}\to Z_{n}\Qd$, \ so the question is reduced from
one about \emph{simplicial} objects over $\C$ to one of algebraic objects,
namely: \ $A_{0}$-algebras. There is a close analagy to the classical theory of
group extensions, where the triviality of an extension \ 
$E:\ 0\to A\to B\to G$ \ is measured by the characteristic class \ 
$\xi(E)\in H^{2}(G;A)$ \ (compare \cite[IV, \S 6]{MacH}). \ 
Similarly, in our case the triviality of the extension is 
measured by the vanishing of a suitable cohomology class in \ 
$H^{n+2}(\Js;\Omega^{n}\Js)$, \ defined as follows\vsm :

Because \ $\jj{n}:\pif Z_{n}\Qd\epic Z_{n}\pif\Qd\cong Z_{n}\Ad$ \ is 
surjective, and \ $\bar{A}_{n+1}$ \ is a free \Pa{\F}, we can choose a 
lifting $\lambda$ in the following diagram:

%
%
%
\begin{figure}[htbp]
\begin{picture}(220,70)(-100,0)
%
%
\put(85,50){$\bar{A}_{n+1}$}
\put(115,55){\vector(1,0){55}}
\put(135,60){$\bar{d}^{A_{n+1}}_{0}$}
\put(175,50){$Z_{n}\Ad$}
\put(205,55){\vector(1,0){25}}
\put(235,50){$0$}
%
%
\multiput(100,45)(0,-3){11}{\circle*{.5}}
\put(100,13){\vector(0,-1){2}}
\put(103,30){$\lambda$}
\put(190,45){\vector(0,-1){30}}
\put(193,30){$\cong$}
%
%
\put(0,0){$0$}
\put(7,5){\vector(1,0){20}}
\put(30,0){$\Omega^{n}\Js$}
\put(60,5){\vector(1,0){23}}
\put(70,10){$i$}
\put(85,0){$\pif Z_{n}\Qd$}
\put(130,5){\vector(1,0){36}}
\put(135,11){$\ju{n}{Q}$}
\put(170,0){$Z_{n}\pif\Qd$}
\put(215,5){\vector(1,0){20}}
\put(240,0){$0$}
\end{picture}
\caption[fig3]{}
\label{fig3}
\end{figure}

\noindent and we can find a map \ $\ell:\bar{Q}_{n+1}\to Z_{n}\Qd$ \ realizing
$\lambda$ \ (again, because \ $\bar{A}_{n+1}=\pif\bar{Q}_{n+1}$ \ is free).
Combined with the ``tautological map'' \ $L_{n+1}\Qd\to M_{n+1}\Qd$ \ (see 
\S \ref{dmat}), which depends only on \ $\tr{n}\Qd$, \ by setting \ 
$Q_{n+1}\DEF \bar{Q}_{n+1}\amalg L_{n+1}\Qd$ \ we obtain \ 
an extension \ $d_{0}:Q_{n+1}\to Q_{n}$ \ of $\ell$ \ (which is a map of \ 
$Q_{0}$-algebras), \ and thus an \ $(n+1)$-truncated simplicial object \ 
$\tr{n+1}\Qd$ \ over $\C$, \ with \ 
$Q_{n+1}\DEF\bar{Q}_{n+1}\amalg L_{n+1}\Qd$, \ and \ 
$\pif\tr{n+1}\Qd\cong \tr{n+1}\Ad$. \ In particular, \ 
$\dd{Q_{n+1}}:C_{n+1}\Qd\to Z_{n}\Qd$ \ induces a map \ $\hat{\lambda}$ \ 
from \ $\pif C_{n+1}\Qd=C_{n+1}\Ad$ \ to \ $\pif Z_{n}\Qd$ \ extending (and 
determined by) the lifting \ $\lambda:\bar{A}_{n+1}\to\pif Z_{n}\Qd$ \ of \ 
$\bar{d}^{A_{n+1}}_{0}$. \ This is a map of \ $A_{0}=\pif Q_{0}$-algebras, by
Assumption \ref{aab}.

Since \ $\ju{n}{Q}\circ(\hat{\lambda}\rest{Z_{n+1}\Ad})=0$, \ the map \ 
$\hat{\lambda}\rest{Z_{n+1}\Ad}$ \ factors through \ 
$\mu:Z_{n+1}\Ad\to\Ker\ju{n}{Q}=\Omega^{n}\Js$, \ and composing $\mu$ with \ 
$\dd{A_{n+2}}:C_{n+2}\Ad\to Z_{n+1}\Ad$ \ defines \ 
$\xi:C_{n+2}\Ad\to \Omega^{n}\Js$ \ -- \ again, a map of \ $A_{0}$-algebras:

%
%
\begin{picture}(360,140)(-60,-15)
%
%
\put(10,90){$C_{n+2}\Ad$}
\put(50,95){\vector(1,0){45}}
\put(63,100){$\dd{A_{n+2}}$}
\put(100,90){$Z_{n+1}\Ad$}
\put(141,95){\vector(1,0){28}}
\put(155,100){$j$}
\put(170,90){$C_{n+1}\Ad$}
\put(210,95){\vector(1,0){42}}
\put(216,100){$\dd{A_{n+1}}$}
\put(255,90){$Z_{n}\Ad$}
\put(285,95){\vector(1,0){25}}
\put(315,90){$0$}
%
%
\put(120,85){\vector(0,-1){73}}
\put(123,50){$\mu$}
\put(180,85){\vector(0,-1){73}}
\put(183,50){$\hat{\lambda}$}
\put(270,85){\vector(0,-1){73}}
\put(275,50){$\cong$}
%
%
\put(40,84){\vector(1,-1){72}}
\put(54,50){$\xi$}
%
%
\put(60,0){$0$}
\put(70,5){\vector(1,0){37}}
\put(110,0){$\Omega^{n}\Js$}
\put(138,5){\vector(1,0){23}}
\put(145,10){$i$}
\put(165,0){$\pif Z_{n}\Qd$}
\put(205,5){\vector(1,0){38}}
\put(215,11){$\ju{n}{Q}$}
\put(250,0){$Z_{n}\pif\Qd$}
\put(290,5){\vector(1,0){20}}
\put(315,0){$0$}
\end{picture}

The cochain \ $\xi=\mu\circ\dd{A_{n+2}}$ \ is clearly a cocycle in the cochain 
complex \ $\Hom_{\RM{\Js}}(\Ad,\Omega\Js)$, \ so it represents a cohomology 
class \ $\chi_{n}\in H^{n+2}(\Js;\Omega^{n}\Js)$, \  called the 
\emph{characteristic class of the extension}.

%
%
\begin{lemma}\label{lthree}\stepcounter{subsection}
The cohomology class \ $\chi_{n}$ \ is independent of the choice of lifting 
$\lambda$.
\end{lemma}

\begin{proof}
Assume that we want to replace $\lambda$ in \S \ref{sco} by a 
different lifting \ $\lambda':\bar{A}_{n+1}\to\pif Z_{n}\Qd$, \ and choose maps \ 
$\ell,\ell':\bar{Q}_{n+1}\to Z_{n}\Qd$ \ realizing \ $\lambda$, \ $\lambda'$ \ 
respectively; \ their extensions to maps \ $Q_{n+1}\to Q_{n}$ \ 
(which we may denote by \ $d_{0}$, \ $d_{0}'$) \ agree on \ $L_{n+1}\Qd$.
We correspondingly having \ 
$\mu':Z_{n+1}\Ad\to\Omega^{n}\Js$ \ and \ $\xi'\DEF\mu'\circ\dd{A_{n+2}}$. \ 

Because \ $Q_{n+1}\DEF \bar{Q}_{n+1}\amalg L_{n+1}\Qd$ \ is a coproduct \ 
of the form \ $\coprod_{i}\M{\alpha_{i}}$, \ by \S \ref{agp} the underlying \ 
group structure on any \ $X\in\C$ \ induces a group structure on \ 
%
\setcounter{equation}{\value{thm}}\stepcounter{subsection}
\begin{equation}\label{eseventeen}
\Hom_{\C}(Q_{n+1},X)
\end{equation}
\setcounter{thm}{\value{equation}}
\noindent (and similarly for \ $\Hom_{\PAlg{\F}}(A_{n+1},\pif X)$).

Therefore, we can set \ $h\DEF (d_{0})^{-1}\cdot(d_{0}'):Q_{n+1}\to Q_{n}$, \ 
and $h$ induces a map \ $\eta:C_{n+1}\Ad\to\pif Z_{n}\Qd$ \ such that \ 
$\eta\rest{\bar{A}_{n+1}}=\lambda^{-1}\cdot \lambda'$. \ Moreover, because \ 
$d_{0}$ \ and \ $d_{0}'$ \ agree outside of \ $\bar{Q}_{n+1}$, \ 
$\ju{n}{Q}\circ\eta=0$. \ Thus \ $\eta$ \ factors through \ 
$\zeta:C_{n+1}\Ad\to\Omega^{n}\Js$, \ which is a map of \ $A_{0}$-algebras 
because \ $\Omega^{n}\Js$ \ is an abelian \ $A_{0}$-algebra (actually, a \ 
$\Js$-module), and $\zeta$ is induced by group operations from the \ 
$A_{0}$-algebra maps \ $\bd$ \ and \ $\bd'$.

Moreover, \ $\zeta\rest{Z_{n+1}\Ad}=\mu-\mu'$ \ in the abelian group structure 
on \ $\Hom_{\RM{\Js}}(-,\Omega^{n}\Js)$ \ (which corresponds to the group
structure of \eqref{eseventeen}).
Thus \ $\xi'-\xi=\hat{\eta}\circ\dd{A_{n+2}}$ \ is a coboundary.
\end{proof}

%
%
\begin{thm}\label{tzero}\stepcounter{subsection}
$\chi_{n}=0$ \ if and only if one can extend \ $\Qd\q{n-1}$ \ to an \ $n$-th
Postnikov approximation \ $\Qd\q{n}$ \ of a resolution of $X$.
\end{thm}

\begin{proof}
First assume that there exists \ $\Yd\simeq \Qd\q{n+1}$ \ with \ 
$\tr{n}\Yd\cong\tr{n}\Qd$: \ by Lemma \ref{ltwo} we know \ 
$\ju{n}{Q}\rest{\Image\ddi{n+1}}$ \ is one-to-one (and onto \ $Z_{n}\pif\Qd$), \ 
for \ $\dd{n+1}:C_{n+1}\Yd\to Z_{n}\Yd=Z_{n}\Qd$, \ and thus \ 
$\Image\ddi{n+1}\cap\Image\partial^{Q}_{n-1}=\{0\}$. \ But then we 
can choose \ $\lambda:\bar{A}_{n+1}\to\pif Z_{n}\Qd$ \ to factor through \ 
$\Image\ddi{n+1}$, \ (and this will induce a map of \ $A_{0}$-algebras because of 
\S \ref{rloop}), so that \ $\mu=0$ \ and thus \ $\xi=0$.

Conversely, if \ $\chi_{n}=0$, \ we can represent it by a coboundary \ 
$\xi=\vartheta\circ\dd{A_{n+2}}$ \ for some \ $A_{0}$-algebra map \ 
$\vartheta:C_{n+1}\Ad\to\Omega^{n}\Js$, \ and thus get \ 
$i\circ\vartheta\rest{\bar{A}_{n+1}}:\bar{A}_{n+1}\to\pif Z_{n}\Qd$. \ 
If we set \ 
$\lambda'\DEF \lambda\cdot (i\circ\vartheta\rest{\bar{A}_{n+1}})^{-1}$, \ 
we have \ $\Image\lambda'\cap\Omega^{n}\Js=\{0\}$. \ We can therefore choose \ 
$\bar{d}_{0}^{Q_{n+1}}:\bar{Q}_{n+1}\to Z_{n}\Qd$ \ realizing \ $\lambda'$, \ 
and then \ $\ddi{Q_{n+1}}$ \ avoids \ 
$\Image(\partial^{Q}_{n-1})\cong\Omega^{n}\Js$, \ so that \ $\tr{n+1}\Qd$ \ 
so constructed yields \ $\Qd\q{n+1}$, \  as required. \ In particular, this 
determines a choice of \ $\Js$-module structure on \ $\Omega^{n+2}\Js$ \ (if
more than one is possible), via \eqref{eeight} \ for \ $n+1$.
\end{proof}

\subsection{notation}
\label{nql}\stepcounter{thm}
If we wish to emphasize the dependence on the choice of $\lambda$, \ we shall
write \ $\Qd\q{n+1}[\lambda]$ \ for the extension of \ $\Qd\q{n}$ \ so 
constructed\vsm .

%
%
\begin{prop}\label{pfour}\stepcounter{subsection}
The class \ $\chi_{n}$ \ depends only on the homotopy type of \ $\Qd\q{n-1}$ \ 
in \ $s\C$.
\end{prop}

\begin{proof}
Assume \ $\Qd\q{n-1}$ \ has been constructed, realizing a simplicial 
resolution of \Pa{\F}s \ $\Ad\to\Js$ \ through simplicial dimension $n$, \ 
and let \ $\Bd\to\Js$ \ be any other \Pa{\F} resolution: we then have a weak 
equivalence \ $\varphi:\Bd\to\Ad$ \ in \ $s\PAlg{\F}$. \ 
Assume by induction on \ $0\leq m<n$ \ that we have constructed an 
$m$-truncated simplicial object \ $\tr{m}\Rd$ \ over $\C$, and a map \ 
$f:\tr{m}\Rd\to\tr{m}\Qd\q{n-1}$ \ realizing \ $\tr{m}\varphi$. \ Moreover,
assume that we have a map of the (split) short exact sequences \eqref{eeight} 
(in dimension $m$) for \ $\Rd$ \ and \ $\Qd$:

%
%
\begin{picture}(220,100)(-70,-15)
%
%
\put(0,50){$0$}
\put(7,55){\vector(1,0){20}}
\put(30,50){$\Omega^{m}\Js$}
\put(60,55){\vector(1,0){23}}
\put(70,60){$i$}
\put(85,50){$\pif Z_{m}\Rd$}
\put(130,55){\vector(1,0){56}}
\put(150,61){$\ju{R}{m}$}
\put(190,50){$Z_{m}\pif\Rd$}
\put(235,55){\vector(1,0){20}}
\put(260,50){$0$}
%
%
\put(40,43){\vector(0,-1){33}}
\put(45,30){$=$}
\put(100,43){\vector(0,-1){33}}
\put(105,30){$(Z_{n}f)_{\#}$}
\put(210,43){\vector(0,-1){33}}
\put(215,30){$Z_{n}(f_{\#})=Z_{n}\varphi$}
%
%
\put(0,0){$0$}
\put(7,5){\vector(1,0){20}}
\put(30,0){$\Omega^{m}\Js$}
\put(60,5){\vector(1,0){23}}
\put(70,10){$i$}
\put(85,0){$\pif Z_{m}\Qd$}
\put(130,5){\vector(1,0){56}}
\put(150,11){$\ju{m}{Q}$}
\put(190,0){$Z_{m}\pif\Qd$}
\put(235,5){\vector(1,0){20}}
\put(260,0){$0$}
\end{picture}

Now, in order to extend $f$ to dimension \ $n+1$, \ we must choose the map \ 
$(\bar{d}_{0}^{R_{m+1}})_{\#}:\pif\bar{R}_{m+1}\to\pif Z_{m}\Rd$ \ (lifting \ 
$\bar{d}_{0}^{B_{m+1}}:\bar{B}_{m+1}\to Z_{m}\Bd$) \ in such a way that \ 
$(Z_{m}f)_{\#}\circ(\bar{d}_{0}^{R_{m+1}})_{\#}=
(\bar{d}_{0}^{R_{m+1}})_{\#}\circ Z_{m}\varphi$. \ Since \ 
$\bar{B}_{m+1}=\pif\bar{R}_{m+1}$ \ is free, it suffices to show that the 
obvious map from \ $\pif Z_{m}\Rd$ \ to the pullback of \ 
$\pif Z_{m}\Qd\xra{\ju{m}{Q}} Z_{m}\pif\Qd=Z_{m}\Ad\xla{Z_{m}\varphi}Z_{m}\Bd$ \ 
is a surjection: \ given \ $(a,b)\in\pif Z_{m}\Qd\times Z_{m}\Bd$ \ with \ 
$\ju{m}{Q}(a)=\varphi(b)$, \ for any \ $z\in\pif Z_{m}\Rd$ \ with \ 
$\ju{m}{R}(z)=b$ \ we have an \ $\omega\in\Omega^{m+1}\Js\subset\pif Z_{m}\Rd$ \ 
such that \ $(Z_{m}f)_{\#}(z\cdot\omega)=(Z_{m}f)_{\#}(z')\cdot\omega=b$ \ in
the diagram above (where \ $\cdot$ \ is the group operation), so \ 
$z\cdot\omega$ \ maps to \ $(a,b)$. \ Thus we can choose \ 
$\bar{d}_{0}^{R_{m+1}}:\bar{R}_{m+1}\to Z_{m}\Rd$ \ in such a way that we can 
define \ $\tr{m+1}\Rd$, \ together with a map \ 
$\tr{m+1}f:\tr{m+1}\Rd\to\tr{m+1}\Qd$ \ realizing \ $\tr{m+1}\varphi$.

Because $\varphi$ was a weak equivalence of 
\emph{resolutions}, it is actually a homotopy equivalence, with homotopy 
inverse \ $\psi:\Ad\to\Bd$, \ say, and the above argument also yields a
homotopy inverse for \ $f\q{m}$  \ (or \ $\tr{m+1}f$). \ Moreover, the
characteristic classes we defined are clearly functorial with respect to
maps in \ $s\C$; \ since the characteristic class \ 
$\chi_{m+1}\in H^{m+3}(\Js;\Omega^{m+1}\Js)$, \ defined for the resolution \ 
$\Ad\to\Js$ \ by means of the lift \ $\bar{d}_{0}^{Q_{m+1}}$, \ must vanish,
by Theorem \ref{tzero}, the same holds for \ $\Rd$, \ so by Theorem \ref{tzero} 
again we can extend \ $\Rd\q{m}$ \ to \ 
$\Rd\q{m+1}$, \ and continue the induction as long as \ $m<n$.
\end{proof}

We deduce the following generalization of Proposition \ref{ptwo}:
%
%
\begin{cor}\label{cone}\stepcounter{subsection}
Given \ $X\in\C$, \ any CW \Pa{\F}\ resolution \ $\Ad\to\pif X$ \ is realizable
as a resolution \ $\Qd\to X$ \ in \ $s\C$.
\end{cor}

One could further extend Proposition \ref{pfour} to obtain a statement about the
naturality of the characteristic classes with respect to morphisms of \Pa{\F}s \ 
$\psi:\Js\to\Ls$. \ However, such a statement would be somewhat convoluted, in
our setting, and it seems better to defer it to a more general discussion of
the realization of simplicial \Pa{\F}s, in \cite{BGoeC}.

\subsection{realization of \Pa{}s}
\label{srpa}\stepcounter{thm}

If \ $G:\Sa\to\G$ \ denotes Kan's simplicial loop functor \ 
(cf.\ \cite[Def.\ 26.3]{MayS}), \ with adjoint \ $\bar{W}:\G\to\Sa$ \ the 
Eilenberg-Mac Lane classifying space functor (cf.\ \cite[\S 21]{MayS}), \ 
and \ $S:\Ta\to\Sa$ \ is the singular set functor, with adjoint \ 
$\|-\|:\Sa\to\Ta$ \ the geometric realization functor 
(see \cite[\S 1,14]{MayS}), then functors
%
\setcounter{equation}{\value{thm}}\stepcounter{subsection}
\begin{equation}\label{esixteen}
\Ta\ \ \substack{S\\ \rightleftharpoons\\ \|-\|}\ \ \Sa\ \ 
\substack{G\\ \rightleftharpoons\\ \bar{W}}\ \ \G
\end{equation}
\setcounter{thm}{\value{equation}}

\noindent induce isomorphisms of the corresponding homotopy categories (see
\cite[I, \S 5]{QuH}), so any homo\-to\-py-theo\-retic question about 
topological spaces may be translated to one in $\G$.
In particular, in order to find a topological space $\X$ having a specified
homotopy \Pa{}\ \ $\Js\cong\pis\X$, \ it suffices to find the corresponding
simplicial group \ $X\in\G$ \ (with the \Pa{\F}\ $\Js$ \ suitably re-indexed).
If \ $\Js$ \ is realizable by such an $X$, any free simplicial resolution \ 
$\Qd\to X$ \ evidently provides a \Pa{}\ resolution \ $\pis\Qd$ \ of \ 
$\Js=\pis X$. \ But the converse is also true: \ if \ $\Qd\in s\G$ \ realizes
some (abstract) \Pa{}\ resolution \ $\Ad\in s\PAlg{}$ \ of \ $\Js$, \ then the
collapse of the Quillen spectral sequence of \cite{QuS}, with 
%
\setcounter{equation}{\value{thm}}\stepcounter{subsection}
\begin{equation}\label{efifteen}
E^{2}_{s,t}=\pi_{s}(\pi_{t}\Qd) \Rightarrow \pi_{s+t}\diag\Qd
\end{equation}
\setcounter{thm}{\value{equation}}

\noindent converging to the diagonal \ $\diag\Qd\in\G$ \ (defined \ 
$(\diag\Qd)_{k}=(Q_{k})^{int}_{k}$) \ implies that \ $\pis\diag\Qd\cong\Js$. \ 
Thus \ $\Js$ \ is realizable by a simplicial group (or topological space) if
and only if some \Pa{}\ resolution \ $\Ad\to\Js$ \ is realizable. 

The characteristic classes \ $(\chi_{n})_{n=0}^{\infty}$ \ 
(whose existence was promised in \cite[\S 1.3]{DKStB} under the name of the 
``$k$-invariants for \ $\Js$), \ thus provide a more succinct (if less 
explicit) version of the theory described in \cite[\S 5-6]{BlaHH} (as 
simplified in \cite[\S 6]{BlaL}), for determining the realizablity of a \Pa{}\ 
in terms of higher homotopy operations \ -- \ which we summarize in

%
%
\begin{thm}\label{tone}\stepcounter{subsection}
Given an (abstract) \Pa{}\ \ $\Js$, \ the following conditions are equivalent:
\begin{enumerate}
\renewcommand{\labelenumi}{(\arabic{enumi})}
\item $\Js$ \ is realizable as \ $\pis\X$ \ for some topological space \ 
$\X\in\Ta$.
\item Any CW \Pa{}\ resolution \ $\Ad\to\Js$ \ is realizable by a simplicial 
space \ $\Qd$.
\item The (inductively defined) characteristic classes \ 
$\chi_{n}\in H^{n+2}(\Js;\Omega^{n}\Js)$ \ ($n=0,1,\dotsc$) \ all vanish.
\end{enumerate}
\end{thm}

Of course, the characteristic class \ $\chi_{n+1}$ \ is determined by the
choice of some extension \ $\Qd\q{n}$ \ of \ $\Qd\q{n-1}$, \ so as usual
our obstruction theory requires back-tracking if at some stage we find \ 
$\chi_{n}\neq 0$. \ We shall now show how we can use other cohomology classes
to determine the choices of extensions at each stage:

\subsection{distinguishing between different resolutions}
\label{sddr}\stepcounter{thm}

A more interesting question, perhaps, is how one can distinguish between
non-equivalent realizations \ $\Qd,\Rd\in s\C$ \ of a fixed
\Pa{\F}\ resolution \ $\Ad\to\Js$ \ of a \emph{realizable} \Pa{\F}\ \ 
$\Js\cong\pif X$. \ Of course, if \ $\Qd$ \ and \ $\Rd$ \ are both 
resolutions (in the resolution model category \ $s\C$) \ of weakly 
equivalent objects \ $X\simeq Y$ \ in the model category \ $\C$, \ then by 
definition \ $\Qd$ \ is weakly equivalent (actually: homotopy equivalent) to \ 
$\Rd$. \ Thus we are looking for a way to distinguish between objects in $\C$,
using the iterative construction of a resolution \ $\Qd\to X$ \ (or 
equivalently, the Postnikov system for \ $\Qd$).

There are a number of possible approaches to this question: one could try to 
construct a homotopy equivalence \ $\Qd\to\Rd$ \ by induction on the Postnikov
tower for \ $\Rd$, \ using an adaptation to \ $s\C$ \ of the classical
obstruction theory for spaces (cf.\ \cite[V, \S 5]{GWhE}). 
Alternatively, one could try directly to construct a map \ $\Qd\to Y$ \ realizing
the augmentation \ $\pif\Ad\to\Js$ \ (see \cite[\S 7]{BlaHH}, and compare 
\cite[\S 5]{BousH}). A description more in this spirit will be given in 
\cite{BGoeC}. 

Here our strategy is similar to that of \S \ref{soer}: rather than assuming that
we are given $X$ and $Y$ to begin with, we try to construct all different 
realizations \ (up to homotopy equivalence in \ $s\C$) \ of a given simplicial
\Pa{\F} \ $\Ad$ \ (which is assumed to be a resolution of a realizable \Pa{\F}\ \ 
$\Js$). \ We start our construction as in \S \ref{sco}, and in the induction 
step we have assume given \ $\tr{n}\Qd$ \ -- \ or equivalently \ 
$\Qd\q{n-1}$, \ satisfying the assumptions of \S \ref{sco} 
(see the proof of Proposition \ref{ptwo}). We ask in how many different ways we
can attach \ $(n+1)$-dimensional ``cells'' to extend the realization one further
dimension.

Again the key lies in the extension of \Pa{\F}s\ of \eqref{eeight}. Of 
course, we may assume that the characteristic class \ 
$\chi_{n}\in H^{n+2}(\Js;\Omega^{n}\Js)$ \ vanishes, so that it is possible
to find ``splittings'' for \eqref{eeight}, given by various liftings 
$\lambda$  in Figure \ref{fig3} \ -- \ all of which yield the same cohomology 
class \ $\chi_{n}$ \ by Lemma \ref{lthree}.
As in the classical case of groups, we find that the difference between two 
such ``semi-direct products'' is represented by suitable cohomology classes, in 
dimension lower by one than the characteristic classes (see 
\cite[IV, \S 2]{MacH}).

\begin{defn}\label{ddoc}\stepcounter{subsection}
Assume given two liftings \ $\lambda,\lambda':\bar{A}_{n+1}\to\pif Z_{n}\Qd$ \ 
in Figure \ref{fig3} above, which define extensions of \ $\tr{n}\Qd$ \ -- \ 
so that, as in the proof of Theorem \ref{tzero}, we may assume without
loss of generality that the corresponding maps \ 
$\mu,\mu':Z_{n+1}\Ad\to\Omega^{n}\Js$ \ vanish. As in the proof of 
Lemma \ref{lthree},  we extend \ $\lambda$, $\lambda'$ \ to face maps \ 
$d_{0},d_{0}':Q_{n+1}\to Q_{n}$, \ define \ $\eta:C_{n+1}\Ad\to\pif Z_{n}\Qd$ \ 
with \ $\ju{n}{Q}\circ\eta=0$, \ and lift to a map of \ $A_{0}$-algebras \ 
$\zeta:C_{n+1}\Ad\to\Omega^{n}\Js$. \ Again \ 
$\zeta\rest{Z_{n+1}\Ad}=\mu-\mu'$, \ which is zero, so \ $\zeta$ is a cocycle 
in \ $\Hom_{\RM{\Js}}(\Ad,\Omega\Js)$, \ representing a cohomology class \ 
$\delta_{\lambda,\lambda'}\in H^{n+1}(\Js,\Omega^{n}\Js)$, \ which we call \ 
the \emph{difference obstruction} for the corresponding Postnikov sections \ 
$\Qd\q{n}[\lambda]$ \ and \ $\Qd\q{n}[\lambda']$ \ (in the notation of 
\S \ref{nql}).
\end{defn}

Just as in the proof of Proposition \ref{pfour}, one can show that the classes \ 
$\delta_{\lambda_{n+1},\lambda'_{n+1}}$ \ in question do not in fact depend on 
the choice of \Pa{\F}\ resolution \ $\Ad\to\Js$, \ but only on the homotopy type 
of \ $\Qd\q{n-1}$ \ in \ $s\C$. \ Their significance is indicated by the following
%
%
\begin{thm}\label{ttwo}\stepcounter{subsection}
If \ $\delta_{\lambda,\lambda'}=0$ \ then the corresponding Postnikov 
sections \ $\Qd\q{n}[\lambda]$ \ and \ $\Qd\q{n}[\lambda']$ \ are weakly 
equivalent.
\end{thm}

\begin{proof}
If \ $\zeta$ \ is a coboundary, there is a map \ 
$\vartheta:C_{n}\Ad\to \Omega^{n}\Js$ \ such that \ 
$\zeta=\vartheta\circ\bd^{A_{n}}$. \ Composing with the inclusion \ 
$i:\Omega^{n}\Js\hra\pif Z_{n}\Qd$ \ yields a morphism of \ $A_{0}$-algebras \ 
$\varphi:A_{n}\to\pif Z_{n}\Qd$. \  If, as in the proof of Proposition 
\ref{ptwo}, we set \ $Q'_{n}\DEF\bar{Q}_{n}\amalg L_{n}\Qd$, \ we may
realize $\varphi$ by a map \ $z':Q'_{n}\to Z_{n}\Qd$. \ Since we assumed \ 
$Q'_{n}$ \ is actually a coproduct of objects in \ $\hF$, \ it is a cogroup 
object in $\C$ by \S \ref{armc}(i), so using the resulting group structure 
on \ $\Hom_{\C}(Q'_{n},Q_{n})$ \ we may set \ 
$s'\DEF k\cdot z:Q'_{n}\to Q_{n}$, \ where \ $k:Q'_{n}\hra Q_{n}$ \ is 
the inclusion. Since $k$ is a trivial cofbration and $Q_{n}$ \ is fibrant in 
$\C$, \ we have a retraction \ $r: Q_{n}\to Q'_{n}$ \ (which is a weak
equivalence). Let \ $s\DEF s'\circ r: Q_{n}\to Q_{n}$. 

Recall from \S \ref{agp} that we have a faithful forgetful functor \ 
$\hat{U}:\C\to\D$, \ where for simplicity we may assume \ $\D=\G$ \ or \ 
$\D=s\RM{R}$ \ (the other cases are trivial). We therefore have a further 
forgetful functor \ $U':\D\to\Ss$, and we denote \ $U'\circ \hat{U}$ \ simply 
by \ $U:\C\to\Ss$. \ The group operation map, while 
not a morphism in \ $\C$ or $\D$, \ is a map \ 
$m:UQ_{n}\times UQ_{n}\to UQ_{n}$ \ in $\Ss$. \ Thus the 
following diagram commutes in $\Ss$:

%
%
\begin{picture}(360,140)(-50,-15)
%
%
\put(0,80){$UZ_{n}\Qd$}
\put(40,85){\vector(1,0){115}}
\put(50,90){$id\,\top\,U(z'\circ r\circ j)$}
\put(160,80){$UZ_{n}\Qd\times UZ_{n}\Qd$}
\put(250,85){\vector(1,0){60}}
\put(255,92){$m\rest{Z_{n}\Qd}$}
\put(315,80){$UZ_{n}\Qd$}
%
%
\put(18,75){\vector(0,-1){33}}
\put(0,58){$Uj$}
\put(205,75){\vector(0,-1){33}}
\put(208,58){$U(j\times j)$}
\put(330,75){\vector(0,-1){30}}
\put(333,58){$Uj$}
%
%
\put(5,30){$UQ_{n}$}
\put(30,35){\vector(1,0){25}}
\put(35,40){$Ur$}
\put(40,25){$\simeq$}
\put(60,30){$UQ'_{n}$}
\put(85,35){\vector(1,0){86}}
\put(90,40){$Uk\,\top\,U(j\circ z')$}
\put(175,30){$UQ_{n}\times UQ_{n}$}
\put(240,35){\vector(1,0){75}}
\put(270,40){$m$}
\put(320,30){$UQ_{n}$}
%
%
\put(25,25){\line(1,-1){15}}
\put(40,10){\line(1,0){265}}
\put(175,-2){$Us$}
\put(305,10){\vector(1,1){15}}
\end{picture}

Since $U$ is faithful, this implies that \ $s\circ j:Z_{n}\Qd\to Q_{n}$ \ 
factors through a map \ $t:Z_{n}\Qd\to Z_{n}\Qd$ \ in $\C$. \ Moreover, 
because we assumed that each \ $\Ma\in\hF$ \ is of the form \ $\Ma=F\Ma'$ \ 
for some \ $\Ma'\in\Ss$ \ (where \ $F=\hat{F}\circ F'$ \ is adjoint to \ 
$U:\C\to\Ss$), \ any map \ $b:\Ma\to Z_{n}\Qd$ \ corresponds under the 
adjunction isomorphism to \ $\hat{b}:\Ma'\to UZ_{n}\Qd$, \ and thus \ 
$t_{\#}\beta=\beta\cdot(\zeta\circ\ju{n}{Q}\beta)$ \ for any \ 
$\beta\in\pif ZS_{n}\Qd$ \ (since the group operation \ $\cdot$ \ in \ 
$\pif Z_{n}\Qd$ \ is induced by $m$ \ -- \ cf.\ \cite[Prop.\ 9.9]{GrayH}).

Now if \ $\ell:\bar{Q}_{n+1}\to Z_{n}\Qd$ \ realizes $\lambda$, \ we have \ 
$(t\circ\ell)_{\#}=(\ell\cdot(z'\circ\ell))_{\#}=
\lambda\cdot(\vartheta\circ\ju{n}{Q}\circ\lambda)=
\lambda\cdot(\lambda^{-1}\cdot\lambda')=
\lambda':\bar{A}_{n+1}\to\pif Z_{n}\Qd$. \ Thus we have a comutative diagram 

%
%
\begin{picture}(360,100)(-70,-16)
%
%
\put(0,50){$\pif\bar{Q}_{n+1}=\bar{A}_{n+1}$}
\put(87,55){\vector(1,0){40}}
\put(100,60){$\lambda'$}
\put(130,50){$\pif Z_{n}\Qd$}
\put(180,55){\vector(1,0){40}}
\put(190,60){$\ju{n}{Q}$}
\put(225,50){$Z_{n}\pif\Qd=Z_{n}\Ad$}
%
%
\put(35,43){\vector(0,-1){30}}
\put(20,28){$id$}
\put(150,43){\vector(0,-1){30}}
\put(153,28){$t_{\#}$}
\put(240,43){\vector(0,-1){30}}
\put(243,28){$id$}
%
%
\put(0,0){$\pif\bar{Q}_{n+1}=\bar{A}_{n+1}$}
\put(87,5){\vector(1,0){40}}
\put(100,10){$\lambda$}
\put(130,0){$\pif Z_{n}\Qd$}
\put(180,5){\vector(1,0){40}}
\put(190,10){$\ju{n}{Q}$}
\put(225,0){$Z_{n}\pif\Qd=Z_{n}\Ad$}
\end{picture}

\noindent which yields a map of \ $(n+1)$-truncated objects \ 
$\rho:\tr{n+1}\Qd[\lambda]\to\tr{n+1}\Qd[\lambda']$ \ (or equivalently, \ 
$\Qd\q{n}[\lambda]\to\Qd\q{n}[\lambda']$). Clearly $\rho$ induces an 
isomorphism in \ $\pi_{k}\pif$ \ for \ $k\leq n+1$. \ 

Now for any choice of lifting $\lambda$ we have \ 
$\pi_{n+2}\pif\Qd\q{n}[\lambda]\cong\Image(\partial^{Q}_{n})$, \ and since \ 
$$
(\vartheta\circ\ju{n}{Q})\rest{\Image(\partial_{n-1}^{Q})}=0,
$$
we find \ $(t_{\#})\rest{\Image(\partial_{n-1}^{Q})}=id$, \ so by \ref{ltwo}(b) 
the diagram

%
%
\begin{picture}(160,95)(-140,-15)
%
%
\put(0,50){$\pi_{\alpha+1}Z_{n}\Qd$}
\put(57,55){\vector(1,0){45}}
\put(70,60){$\partial^{Q}_{n}$}
\put(105,50){$\pia Z_{n+1}\Qd$}
%
%
\put(25,42){\vector(0,-1){33}}
\put(10,28){$t_{\#}$}
\put(130,42){\vector(0,-1){33}}
\put(133,28){$id$}
%
%
\put(0,0){$\pi_{\alpha+1}Z_{n}\Qd$}
\put(57,5){\vector(1,0){45}}
\put(70,10){$\partial^{Q}_{n}$}
\put(105,0){$\pia Z_{n+1}\Qd$}
\end{picture}

\noindent commutes. Thus $\rho$ induces an isomorphism on \ 
$\Image(\partial^{Q}_{n})$, \ so that \ 
$(\rho)_{\ast}:\Qd\q{n}[\lambda]\to\Qd\q{n}[\lambda']$ \ is a weak 
equivalence.
\end{proof}

\begin{remark}\label{rkin}\stepcounter{subsection}
Given a (realizable) \Pa{\F}\ \ $\Js$, \ a CW resolution \ 
$\Ad\in s\PAlg{\F}$ \ of \ $\Js$, \ and a fixed (but arbitrary) choice 
object \ $X\in\C$ \ with \ $\pif X\cong\Js$, \ by Corollary \ref{cone} we have
a corresponding resolution \ $\Qd\to X$. \ If \ $X'\in\C$ \ is another 
realization of \ $\Js$ \ with corresponding \ $\Qpd\to X'$, \ we may 
assume without loss of generality that \ 
$\Yd'\DEF(\Qpd)\q{n}\simeq\Yd\DEF\Qd\q{n}$ \ for some \ $n\geq 0$, \ 
with \ $\lambda,\lambda':A_{n+2}\to\pif Z_{n+1}\Qd\cong\pif Z_{n+1}\Qpd$ \ the 
respective liftings. 
\end{remark}

\subsection{different realizations of a \Pa{}}
\label{sdrp}\stepcounter{thm}

Assume given an abstract \Pa{} \ $\Js$, \ which is known to be realizable
(e.g., by the cohomological criterion of Theorem \ref{tone}). We wish
to distinguish between the various non-weakly equivalent
realizations of \ $\Js$ \ by topological spaces (or simplicial groups). 
The spectral sequence \eqref{efifteen} implies that in order for two
such \ $X,X'\in\G$ \ (with \ $\pis X\cong\Js\cong \pis X'$) \ to be weakly
equivalent, it suffices that their corresponding resolutions \ $\Qd\to X$ \ 
and \ $\Qpd\to X'$ \ be weakly equivalent (and thus homotopy equivalent) in the
resolution model category.   This is in fact the main reason for considering
this model category structure on \ $s\G$ \ in the first case (and justifies
its original name of \ ``$E^{2}$-model category'' in \cite{DKStE}).

Note, however, that this is not a necessary condition;  an alternative
model structure on \ $s\Ss$ \ (or \ $s\G$), \ defined in \cite{MoerB}, has as 
weak equivalences precisely those maps in \ $s\C$ \ inducing an equivalence
on the realizations.

The difference obstructions \ $\delta_{\lambda,\lambda'}$, \ which yield an
inductive procedure for distinguishing between various realizations of 
a given \Pa{}\ resolution \ $\Ad\to\Js$, \ thus again provide an alternative
to the theory described in \cite[\S 7]{BlaHH} (as simplified in 
\cite[\S 4.9]{BlaHO}) for distinguishing between different realizations
of a given \Pa{}, in terms of higher homotopy operations.

To state this explicitly, assume given an (abstract) \Pa{}\ \ $\Js$, \ 
a CW resolution \ $\Ad\in s\PAlg{}$ \ of \ $\Js$, \ and two realizations \ 
$\Qd,\Qpd\in s\G$ \ of \ $\Ad$, \ determined as in \S \ref{sddr} by 
successive choices of lifts \ $\lambda_{k+1}:\bar{A}_{k+1}\to\pif Z_{k}\Qd$ \ 
and \ $\lambda'_{k+1}:\bar{A}_{k+1}\to\pif Z_{k}\Qpd$. \ By \S \ref{srpa},
we know that the realizations \ $X\DEF\diag \Qd$ \ and \ $X'\DEF\diag \Qpd$ \ 
are two realizations of \ $\Js$. \ If \ 
$\delta_{\lambda_{0},\lambda'_{0}}=0$, \ there is a weak equivalence \ 
$f_{0}:(\Qpd)\q{0}\simeq\Qd\q{0}$, \ which we can use to push forward \ 
$\lambda'_{1}:\bar{A}_{2}\to\pif Z_{1}\Qpd$ \ to \ 
$\lambda''_{1}:\bar{A}_{2}\to\pif Z_{1}\Qd$ \ so it is meaningful to consider
$\delta_{\lambda_{1},\lambda'_{1}}\DEF\delta_{\lambda_{1},\lambda''_{1}}
\in H^{2}(\Js,\Omega\Js)$. \ Proceeding in this way we obtain the following

%
%
\begin{thm}\label{tthree}\stepcounter{subsection}
Assume given a \Pa{}\ \ $\Js$, \ a CW resolution \ 
$\Ad\in s\PAlg{}$ \ of \ $\Js$, \ and two topological spaces \ 
$\X,\X'\in\Ta$ \ realizing \ $\Js$, \ corresponding to \ $X,X'\in \G$ \ under 
\eqref{esixteen}. \ Let \ $\Qd,\Qd'\in s\G$ \ be 
CW resolutions of \ $X,X'$ \ respectively, \ determined as in \S 
\ref{sddr} by successive choices of lifts \ 
$\lambda_{n+1}:\bar{A}_{n+1}\to\pif Z_{n}\Qd$ \ and \ 
$\lambda'_{n+1}:\bar{A}_{n+1}\to\pif Z_{n}\Qd'$. \ If the 
difference obstructions \ 
$\delta_{\lambda_{n+1},\lambda'_{n+1}}\in H^{n+2}(\Js,\Omega^{n+1}\Js)$ \ 
vanish for all \ $n\geq 0$, \ then \ $\X$ \ and \ $\X'$ \ are weakly 
equivalent.
\end{thm}

Again, these classes satsify certain naturality conditions, which are more easily
stated for simplicial \Pa{\F}s: see \cite{BGoeC}.

\begin{remark}\label{rdpr}\stepcounter{subsection}
Theorem \ref{tthree} provides a collection of algebraic invariants \ -- \ 
starting with the homotopy \Pa{}\ \ $\pis\X$  \ -- \ for distinguishing between
(weak) homotopy types of spaces. As with the ordinary Postnikov systems and
their $k$-invariants, these are not actually invariant, in the sense that
distinct values (i.e., non-vanishing difference obstructions) do not guarantee
distinct homotopy types. Thus we are still far from a full algebraization of
homotopy theory \ -- \ even if we disregard the fact that \Pa{}s, not too 
mention their cohomology groups, are rather mysterious objects, and no 
non-trivial naturally occurring examples are fully known to date. 

Note, however, that we have a considerable simplification of the theory in the
case of the rational homotopy type of simply-connected spaces: in this
case the \Pa{\F}s in question are just connected graded Lie algebras over 
$\Q$, and the cohomology theory reduces to the usual 
Cartan-Eilenberg cohomology of Lie algebras. The obstruction theory we define 
appears to be the Lie algebra version of the theory for graded algebras
due to Halperin and Stasheff in \cite{HStaO}. See also \cite[\S III]{OukiH} \ 
and \ \cite{FelDT}.

Another such simplification occurs when we consider only the stable homotopy 
type: in this case \Pa{\F}s are just graded 
modules over the stable homotopy ring \ ${\mathbb \pi}\DEF\pis^{S}S^{0}$, \ 
and the cohomology groups in question are \ 
$\Ext^{\ast}_{\mathbb \pi}(\Js,\Sigma^{-n}\Js)$. \ Here we have no action of the
fundamental group to worry about.

Furthermore, the spectral sequence of \eqref{efifteen} implies that if \ 
$\Qd\q{n}\cong(\Qpd)\q{n}$, \ then \ also \ 
$(\diag\Qd)\q{n}\cong(\diag\Qpd)\q{n}$, \ so one can also use the theory
described above ``within a range''.
\end{remark}

\end{document}